\documentclass{autart}               

\usepackage[round]{natbib}

\usepackage{graphicx}          
\usepackage{savesym}
\savesymbol{AND}
\usepackage{algorithmic }
\usepackage{amsmath,amssymb}
\usepackage{bm}
\usepackage{url}
\usepackage{multirow}
\usepackage{color}
\usepackage[normalem]{ulem}
\usepackage{cancel}

\newcommand{\cut}[1]{}
\newcommand{\pc}[1]{{\color{blue}#1}}
\newcommand{\zs}[1]{{\color{red}#1}}

\newtheorem{asu}{Assumption}
\newtheorem{lemma}{Lemma}
\newtheorem{remark}{Remark}
\newtheorem{definition}{Definition}
\newtheorem{theorem}{Theorem}
\newtheorem{corollary}{Corollary}
\allowdisplaybreaks

\begin{document}

\begin{frontmatter}
%
\title{On convergence rates of game theoretic reinforcement learning algorithms} 
%
\thanks[footnoteinfo]{Z. Hu and M. Zhu were partially supported by  ARO W911NF-13-1-0421 (MURI), NSA H98230-15-1-0289 and NSF ECCS-1710859. P. Chen and P. Liu were partially supported by ARO W911NF-13-1-0421 (MURI) and NSF CNS-1422594.}
%
\author[EE]{Zhisheng Hu}\ead{zxh128@psu.edu},    
\author[EE]{Minghui Zhu}\ead{muz16@psu.edu},               
\author[JD]{Ping Chen}\ead{ pzc10@ist.psu.edu},  
\author[IST]{Peng Liu}\ead{pliu@ist.psu.edu}

\address[EE]{School of Electrical Engineering and Computer Science, Pennsylvania State University, 201 Old Main, University Park, PA, 16802, USA}  
\address[JD]{JD.com, No. 18 Kechuang 11 Street, BDA, Beijing,10111, China} 
\address[IST]{College of Information Sciences and Technology, Pennsylvania State University, 201 Old Main, University Park, PA, 16802, USA}             

\begin{keyword}                           
         Distributed control;
       Game theory;
        Learning in games
\end{keyword}                             

\begin{abstract}                          
This paper investigates a class of multi-player discrete games where each player aims to maximize its own utility function. Each player does not know the other players' action sets, their deployed actions or the structures of its own or the others' utility functions. Instead, each player only knows its own deployed actions and its received utility values in recent history. \cut{Second, utility values are subject to I.I.D. random noises.} We propose a reinforcement learning algorithm which converges to the set of action profiles which have maximal stochastic potential with probability one. Furthermore, the convergence rate of the proposed algorithm is quantified.\cut{ When the interactions of the players consist of a weakly acyclic game and the received utilities are subject to uniformly bounded random noises, the convergence to the set of \pc{anti-$\mathcal{E}$-Nash equilibria} is guaranteed.} The algorithm performance is verified using two case studies in the smart grid and cybersecurity.
\end{abstract}

\end{frontmatter}

\section{Introduction}
Game theory provides a mathematically rigorous framework for multiple players to reason about each other. In recent years, game theoretic learning has been increasingly used to control large-scale networked systems due to its inherent distributed nature. In particular, the network-wide objective of interest is encoded as a game whose Nash equilibria correspond to desired network-wide configurations. Numerical algorithms are then synthesized for the players to identify Nash equilibria via repeated interactions. Multi-player games can be categorized into discrete games and continuous games. In a discrete (resp. continuous) game, each player has a finite (resp. an infinite) number of action candidates. As for discrete games, learning algorithms include best-response dynamics, better-response dynamics, factitious play, regret matching, logit-based dynamics and replicator dynamics. Please refer to \citep{TB-GO:99,fudenberg1998theory,sandholm2010population,young2001individual} for detailed discussion. As an important class of continuous games, generalized Nash games were first formulated in \cite{Arrow.Debreu:54}, and see survey paper \cite{Facchinei2007} for a comprehensive exposition. A number of algorithms have been proposed to compute generalized Nash equilibria, including, to name a few, ODE-based methods \cite{Rosen:65}, nonlinear Gauss-Seidel-type approaches \cite{Pang.Scutari.Facchinei.Wang:08}, iterative primal-dual Tikhonov schemes \cite{Yin.Shanbhag.Mehta:11}, and best-response dynamics \cite{Palomar.Eldar:10}. Game theory and its learning have found many applications; e.g., traffic routing in Internet~\cut{\citep{MHM-QZ-TA-TB:13,1008358}}\cite{1008358}, urban transportation~\cite{Roumboutsos2008209}, mobile robot coordination~\citep{GA-JRM-JSS:07a\cut{,GT-JB:08},hatanaka2016payoff} and power markets~\citep{6425983,Zhu:PESGM14}.\cut{\zs{ (to be finished)}}

In many applications, players can only access limited information about the game of interest. For example, each player may not know the structure of its own utility function. Additionally, during repeated interactions, each player may not be aware of the actions of other players. These informational constraints motivate recent study on  payoff-based or reinforcement learning algorithms where the players adjust their actions only based on their own previous actions and utility measurements. The papers~\citep{doi:10.1137/070680199,MZ-SM:SIAM13,hatanaka2016payoff} study discrete games, and their approaches are based on stochastic stability~\cite{FOSTER1990219}. As mentioned in Remark 3.2 of~\cite{MZ-SM:SIAM13}, paper~\cite{doi:10.1137/070680199} proposes an algorithm to find Nash equilibrium of weakly acyclic games with an arbitrarily high probability by choosing an arbitrarily small and fixed exploration rate in advance. The analysis in~\cite{doi:10.1137/070680199} is based on homogeneous Markov chains and more specifically the theory of resistance trees~\cite{young1993evolution}. \cite{MZ-SM:SIAM13} extends the results in~\cite{doi:10.1137/070680199} by adopting diminishing exploration rates and ensures convergence to Nash equilibrium and global optima with probability one. The analysis of~\cite{MZ-SM:SIAM13} is based on based on strong ergodicity of inhomogeneous Markov chains. As for continuous games, the papers~\citep{PF-MK-TB:11,SL-MK:11,MS-KJ-DS:12} employ extremum seeking\cut{ and the paper~\cite{QZ-HT-TB:13} uses finite-difference approximations to estimate unknown gradients} and the paper~\cite{MZ-EF:AUTOMATICA16} uses finite-difference approximations to estimate unknown partial (sub)gradients. Notice that all the aforementioned papers focus on asymptotic convergence and none of them quantifies convergence rates.

\emph{Contribution}: In this paper, we study a class of multi-player discrete games\cut{. The particular challenges we consider include (i)} where each player is unaware of\cut{ the actions taken by the others and its own utility function} the other players' action sets, their deployed actions or the structures of its own or the others' utility functions. \cut{(ii) its received utilities are subject to independent and identically distributed (I.I.D.) random noises.} \cut{Inspired by~\cite{MZ-SM:SIAM13}, }We propose a reinforcement learning algorithm where, at each iteration, each player, on one hand, exploits successful actions in recent history via comparing received utility values, and on the other hand, randomly explores any feasible action with a certain exploration rate.
The algorithm is proven to be convergent to the set of action profiles with maximum stochastic potential with probability one. Furthermore, an upper bound on the convergence rate is derived and  is minimized when the exploration rates are restricted to $\mathbf{p}$-series.\cut{provides a guideline to choose the exploration rates to speed up the algorithm.}\cut{Furthermore, the convergence rate is formally quantified.} When the interactions of the players consist of a weakly acyclic game and the received utilities are subject to uniformly bounded random noises, the convergence to the set of pure anti-$\mathcal{E}$-Nash equilibria is guaranteed. The algorithm performance is verified using two case studies in the smart grid and cybersecurity.\cut{ A set of numerical simulations based on two real-world applications are used to evaluate our algorithm.} A preliminary version of this paper was published in~\cite{zhu2014reinforcement} where convergence rates and measurement noises are not discussed. Further,~\cite{zhu2014reinforcement} focuses on the application on adaptive cyber defense, and this paper focuses on theory of learning in games. The analysis of two papers is significantly different.

\cut{Hence in this paper, we do not specify the strategy space of the attacker or the utility function. Instead, we propose a defense algorithm that enables the defender to update its defense action based on the feedbacks online.}

\section{Problem formulation and learning algorithm}\label{model}

In this section, we introduce a class of multi-player games where the information each player accesses is limited. Then, we present a learning algorithm under which the action profiles of the players converge to the set of action profiles which have maximum stochastic potential.
\subsection{Game formulation}\label{general_model}
\begin{figure}
  \centering
  \includegraphics[width = \linewidth]{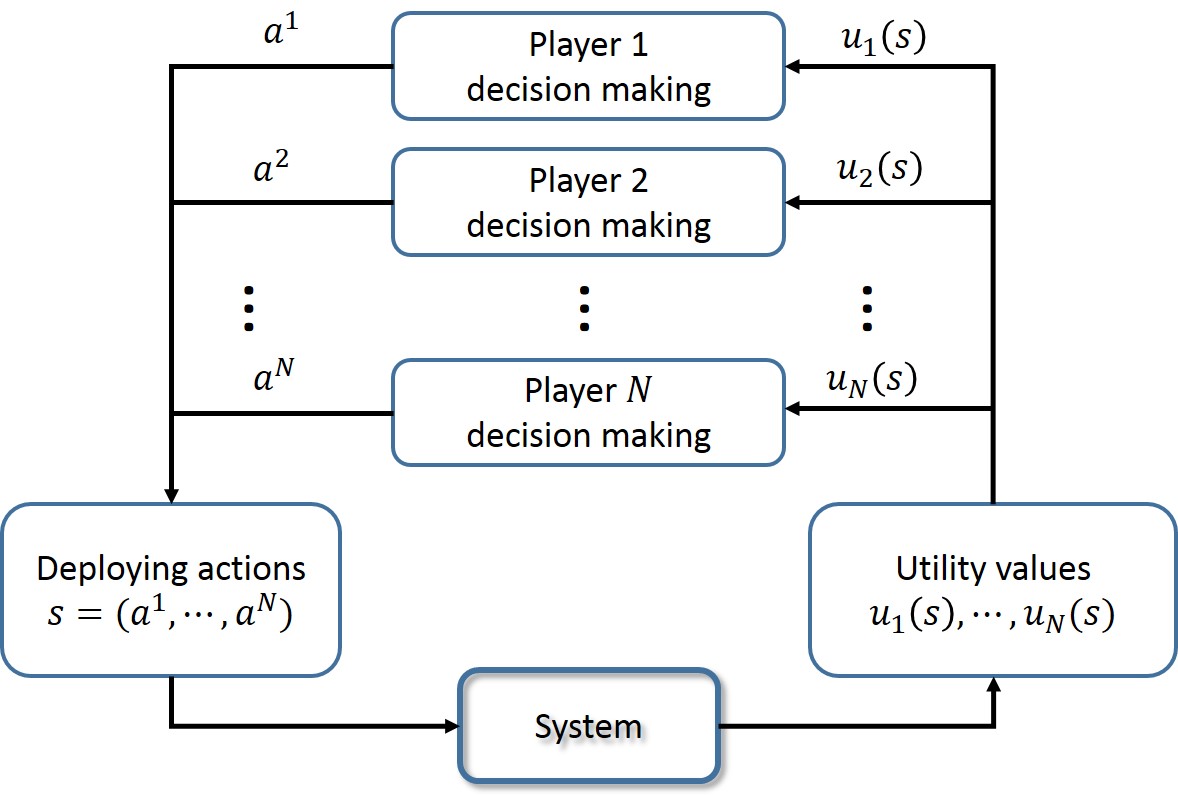}
  \caption{Game model}\label{fig1}

\end{figure}

The system model in Figure~\ref{fig1} characterizes the interactions of $N$ players in a non-cooperative game. Each component in the figure will be discussed in the following paragraphs.

{\bf Players.}
We consider $N$ players $\mathcal{V}\triangleq\{1,\cdots,N\}$ and each player has a finite set of actions. Let $\mathcal{A}_i$ denote the action set of player $i$ and $a^i\in \mathcal{A}_i$ denote an action of player $i$. Denote $\mathcal{S}\triangleq \mathcal{A}_1\times\cdots\times\mathcal{A}_N$ as the Cartesian product of the action sets, where $s\triangleq (a^1,\cdots,a^N)\in \mathcal{S}$ is denoted as an action profile of the players.

{\bf Utility.}
Under the influence of an action profile, the system generates a utility value for each player. The utility function for player $i\in\mathcal{V}$ is defined as $u_i:\mathcal{S}\to \mathbb{R}$.  At the end of iteration $t$, the utility value $u_i(t)=u_i(s(t))$ is measured and sent to player $i$. The utility value received by player $i$ is denoted by $\tilde u_i (t)\triangleq u_i(t)+w_i(t)$, where $w_i(t)$ is the measurement noise.

\cut{{\bf Synchronous update.}
In this paper, we assume that each player periodically and synchronously updates its action based on its feedback. In particular, each player performs mate-analysis and then deploys an action until next update. The time interval between two updates is denoted as a iteration; i.e., time instance $t$ represents the $t^{th}$ iteration. We assume that the lengths of the iterations are uniform.}

{\bf Informational constraint.}
Each player does not know the other players' action sets or their deployed actions. Besides, each player is unaware of the structure of its own or the others' utility functions. At iteration $t$, each player only knows its deployed actions and its received utility values in the past; i.e., $a^i(0), \cdots, a^i(t-1), \tilde u_i(0), \cdots, \tilde u_i(t-1)$.

The above informational constraint has been studied in several recent papers. For example, the authors in \citep{MZ-SM:SIAM13,MS-KJ-DS:12,hatanaka2016payoff} investigate coverage optimization problems for mobile sensor networks where mobile sensors are unaware of environmental distribution functions.\cut{ The authors in \cite{PF-MK-TB:11} study the non-cooperative game in duopoly market structure where both the players only need to measure their own payoff values.} The authors in \cite{6515297} study the problem of optimizing energy production in wind farms where each turbine knows neither the functional form of the power generated by the wind farm nor the choices of other turbines. The authors in \citep{PF-MK-TB:11,MZ-EF:AUTOMATICA16} consider convex games where each player cannot access its game components.\cut{ In addition, the defender knows the class of the zero-day attack; e.g., buffer over-read or code reuse or combination, and then it can avoid employing ill-matching defense techniques.}

\subsection{Problem statement}%
Under the above informational constraint, we aim to synthesize a learning algorithm under which the action profiles of the players converge to the set of action profiles with maximum stochastic potential. \cut{In contrast to existing work, our algorithm considers heterogeneous exploration rates and analyzes convergence property when measurements are subject to independent and identically distributed (I.I.D.) random noises.} We will quantify the convergence rate of the proposed algorithm in contrast to asymptotic convergence in existing work.

\subsection{Learning algorithm}
Inspired by~\cite{MZ-SM:SIAM13}, we propose a learning algorithm called the RL algorithm, where each player updates its actions only based on its previous actions and its received utility values. On the one hand, each player chooses the most successful action in recent history. It represents the exploitation phase.\cut{On one hand, the defender exploits its previous actions and chooses the most successful action in recent history.} However, the exploitation is not sufficient to guarantee that the player can choose the best action given others'.\cut{choose the most successful action within its whole action set.} So on the other hand, the player uniformly chooses one action from its action set. It represents the exploration phase. The specific update rule is stated in the RL algorithm. At iterations $t=0$ and $t=1$, each player uniformly chooses one action from its action set as initialization (Line 3). Starting from iteration $t=2$, with probability $1-\tilde\epsilon_i(t)$,  player $i$ chooses the action which generates a higher utility value in last two iterations as current action (Line 8-13). This represents the exploitation where player $i$ reinforces its previous successful actions. With probability $\tilde\epsilon_i(t)$, player $i$ uniformly selects an action from its action set $\mathcal{A}_i$ (Line 14). This represents the exploration and makes sure that each action profile is selected infinitely often.
Note that ${\rm sample}(\mathcal{A}_i)$ in Line 14 represents uniformly choosing one element from set $\mathcal{A}_i$.



\begin{algorithm}[h]\caption{Reinforcement learning (RL) algorithm} \label{algo:general}
\begin{algorithmic}[1]
\WHILE {$0\le t \le 1$}{
\FOR {$i\in\mathcal{V}$}\STATE{
$a^i(t) \leftarrow {\rm sample}(\mathcal{A}_i)$;
}\ENDFOR
}
\ENDWHILE
\WHILE{$t\ge 2$}{
\FOR {$i\in\mathcal{V}$}{
\STATE With prob. $(1-\tilde\epsilon_i(t))$, \IF{$\tilde u_i(t-1)\ge\tilde u_i(t-2)$}\STATE{$a^i(t) = a^i(t-1)$;}
\ELSE\STATE{$a^i(t) = a^i(t-2)$;}\ENDIF

\STATE With prob. $\tilde\epsilon_i(t)$, $a^i(t) \leftarrow{\rm sample}(\mathcal{A}_i\cut{\setminus \{d(t-1),d(t-2)\}})$;
}\ENDFOR
}\ENDWHILE
\end{algorithmic}
\end{algorithm}

\cut{Note that in the RL algorithm, we include the algorithm of the attacker to make the interactions complete. We do not specify the exploitation algorithm $ALG_a$ of the attacker.\cut{ Instead, we simply assume that the attacker exploits its previous actions and randomly explores the attack action set $\mathcal{A}$ with a certain probability.}}

\section{Analysis}\label{sec:error}
In this section, we will present the analytical results of the RL algorithm.\cut{ The results focus on the convergence properties of the defense and attack actions: (1) whether they converge, (2) where they converge, and (3) how fast they converge.}
\subsection{Notations and assumptions}\label{subsec:notation}
We first introduce the notations and assumptions used throughout the paper. Denote by $|\mathcal{V}|$ the cardinality of player set, $ |\mathcal{A}_i|$ the cardinality of action set of player $i$ and $|\mathcal{A}|_{\infty}\triangleq \max\limits_{i\in \mathcal{V}}|\mathcal{A}_i|$ the maximum cardinality among all action sets.
The exploration rate for player $i$ at iteration $t$ is decomposed into two parts; i.e., $\tilde\epsilon_i(t)\triangleq\epsilon_i(t)+{e_i}(t)$, where $\epsilon_i(t)=\gamma_i{\epsilon}^c(t)$, $\gamma_i>0$, ${\epsilon}^c(t)$ is common for all the players and ${e_i}(t)$ represents the exploration deviation.
Define $e(t)\triangleq({e_1}(t),\cdots,{e_N}(t))^T$, $\tilde{\epsilon}(t) \triangleq (\tilde\epsilon_1(t),\cdots, \tilde\epsilon_N(t))^T$ and ${\epsilon}(t) \triangleq (\epsilon_1(t),\cdots, \epsilon_N(t))^T$. And we define $e_r(t)\triangleq||{e}(t)||_\infty^N/\prod\limits_{i=1}^N\tilde\epsilon_i(t)$. Here we denote by $||\cdot||_\infty$ the infinity norm of a vector. In addition, we also use $||\cdot||$ to represent the $L^1$-norm of a vector, and $||P|| \cut{\triangleq\max\limits_{k}\bm{\sum}\limits_i|p_{ik}|}$ to represent the 1-norm of a matrix. We assume the measurement noises satisfy:

\begin{asu}\label{asu_error_o}
For each $i\in\mathcal{V}$, $\{w_i(t)\}_{t\ge 0}$ is a sequence of real-valued random variables that are independent and identically distributed (I.I.D.).


\end{asu}


I.I.D. random noises are widely adopted in distributed control and optimization problems and see~\citep{XIAO200733,doi:10.1137/06067359X}. In addition, we assume the exploration rates satisfy:

\begin{asu}\label{asu6}
(1). For each $i\in\mathcal{V}$, $\epsilon_i(t)\in(0,1]$ is non-negative, strictly decreasing, and $\lim \limits_{t\to \infty}\epsilon_i(t)=0$.
(2). For each $i\in\mathcal{V}$, the sequences $\{\prod\limits_{i=1}^N\epsilon_i(t)\}$ and $\{\prod\limits_{i=1}^N\tilde\epsilon_i(t)\}$ are not summable. (3). $\lim\limits_{t\to\infty}e_r(t)=0$.
\end{asu}

Assumption~\ref{asu6} indicates that the players can choose heterogeneous exploration rates. The exploration rates diminish slowly enough and their deviations decrease in  faster 
rates than the common part.
In the paper~\cite{MZ-SM:SIAM13\cut{doi:10.1137/070680199, hatanaka2016payoff}}, it is assumed that exploration rates $\epsilon_i(t)$ are identical for all $i$, diminishing and not summable. Assumption~\ref{asu6} allows for heterogeneous exploration rates and includes homogeneous exploration rates in the paper~\cite{MZ-SM:SIAM13\cut{doi:10.1137/070680199,hatanaka2016payoff}} as a special case. Actually, papers \citep{\cut{doi:10.1137/090770102,}6286992,6580575} adopt heterogenous step-sizes for distributed optimization and game theory. They impose similar assumptions on the step-sizes.\cut{ Assumption~\ref{asu6} indicates that as the game goes on, all the players' exploration rates decrease in the same order to 0 and the exploration rates do not decrease too fast.}\cut{ This assumption is inspired by simulated annealing~\cite{aarts1988simulated}.}

{\bf Markov chain induced by the RL algorithm.} Denote by $\mathcal{Z}\triangleq\mathcal{S}\times\mathcal{S}$ the state space, where each state $z(t)\triangleq(s(t),s(t+1))$ consists of the action profiles at iteration $t$ and the next iteration. And denote by $diag(\mathcal{S}\times\mathcal{S})\triangleq\{(s,s)|s\in \mathcal{S}\}$ the diagonal space of $\mathcal{Z}$. By the definition of $z(t)$, the sequence $\{z(t)\}_{t\ge 0}$ forms a time-inhomogeneous Markov chain, denoted by $\mathcal{M}$. We define $P^{\tilde{\epsilon}(t)}$ as the transition matrix of Markov chain $\mathcal{M}$ at iteration $t$, where each entry $P^{\tilde{\epsilon}(t)}(z',z)$ represents the transition probability from state $z'$ to $z$. Besides, denote by $\pi(t)$ the distribution on $\mathcal{Z}$ at iteration $t$. 

{\bf $z$-tree of time-homogenous Markov chain $\mathcal{M}^{\tilde{\epsilon}}$.} Given any two distinct states $z'$ and $z$ of Markov chain $\mathcal{M}^{\tilde{\epsilon}}$, consider all paths starting from $z'$ and ending at $z$. Denote by $p_{z'z} $ the largest probability among all possible paths from $z'$ to $z$. A path might contain\cut{ some $i\ge 0$} intermediate states $z_{1},\cdots,z_{k}$ ($k=0$ means there is no intermediate state) between $z'$ and $z$. So $p_{z'z}$ is the product of $P^{\tilde{\epsilon}}(z',z_1),P^{\tilde{\epsilon}}(z_1,z_2),\cdots, P^{\tilde{\epsilon}}(z_k,z)$. We define graph $\mathcal{G}(\tilde{\epsilon})$ where each vertex of $\mathcal{G}(\tilde{\epsilon})$ is a state $z$ of Markov chain $\mathcal{M}^{\tilde{\epsilon}}$ and the probability on edge $(z',z)$ is $p_{z'z}$. A $z$-\emph{tree} on $\mathcal{G}(\tilde{\epsilon})$ is a spanning tree rooted at $z$ such that from every vertex $z'\neq z$, there is a unique path from $z'$ to $z$. Denote by $G_{\tilde{\epsilon}}(z)$ the set of all $z$-\emph{trees} on $\mathcal{G}(\tilde{\epsilon})$. The total probability of a $z$-\emph{tree} is the product of the probabilities of its edges. The \emph{stochastic potential} of the state $z$ is the largest total probability among all $z$-\emph{trees} in $G_{\tilde{\epsilon}}(z)$. Let $\Lambda(\tilde{\epsilon})$ be the states which have maximum stochastic potential for a particular $\tilde{\epsilon}\in(0,1]$. Denote the limit set $\Lambda^\ast\triangleq \lim\limits_{\tilde{\epsilon}\to 0}\Lambda(\tilde{\epsilon})$. And the elements in $\Lambda^\ast$ are referred to as \textit{stochastically stable states}.\cut{ Any of such trees is denoted by $T_{\max}(\tilde{\epsilon}(t_q))\in G_q(z)$.}

\cut{
{\bf $z$-tree of Markov chain $\mathcal{M}_q$ at defense iteration $t_q$.} For any iteration $t_q$ and two distinct states $z',z\in \mathcal{Z}$ of the corresponding Markov chain $\mathcal{M}_q$, denote $p^q_{z'z} $ as the largest probability among all possible paths from $z'$ to $z$. A path from $z'$ to $z$ might contain some $i\ge 0$ middle states $z_{1},\cdots,z_{i}$ ($i=0$ means there is no middle state). So $p^q_{z'z}$ is the product of $P_q^{\tilde{\epsilon}(t_q)}(z',z_1),P_q^{\tilde{\epsilon}(t_q)}(z_1,z_2),\cdots, P_q^{\tilde{\epsilon}(t_q)}(z_i,z)$. Define a set of complete directed graphs $\{\mathcal{G}(t)\}$. At $t=t_q$, each vertex of $\mathcal{G}(t_q)$ is a state $z$ of Markov chain $\mathcal{M}_q$ and the probability on edge $(z',z)$ is $p^q_{z'z}$. Define a $z$-\emph{tree} on $\mathcal{G}(t_q)$ as a spanning tree such that from every vertex $z'\neq z$, there is a unique path from $z'$ to $z$. Denoted by $G_q(z)$ the set of all $z$-\emph{trees} on $\mathcal{G}(t_q)$. The total probability of a $z$-\emph{tree} is the product of the probabilities of its edges. The \emph{stochastic potential} of the state $z$ is the large total probability among all $z$-\emph{trees} in $G_q(z)$ . Let $\Lambda(\tilde{\epsilon}(t))$ be the set of states which have maximum stochastic potential at particular defense iteration $t$. And  denote $\Lambda^\ast$ the set of states which have maximum stochastic potential as $t\to \infty$.\cut{ Any of such trees is denoted by $T_{\max}(\tilde{\epsilon}(t_q))\in G_q(z)$.}
}
\begin{remark}
The above notions are inspired by the resistance trees theory~\cite{young1993evolution}. However, the above notions are defined for any $\tilde{\epsilon}\in (0,1]$ instead of $\tilde{\epsilon}\to 0$ in the resistance trees theory. This allows us to characterize the transient  performance of the RL algorithm. \qed
\end{remark}

\cut{{\bf Stationary distribution.} Denote by $\pi(t)$ the distribution on $\mathcal{Z}$ at defense iteration $t$\cut{ and $\pi(0)$ the initial distribution on $\mathcal{Z}$ of Markov chain $\{P_t\}$}. Define stochastic vector $\pi^\ast(t_q)$ as the stationary distribution of the Markov chain $\mathcal{M}_q$ when $t=t_q$. And by Theorem 3.1 in Chapter 6 of \cite{freidlin2012random}, $\pi^\ast(t_q)$ can be represented as follows.
\begin{align}
\pi^\ast(t_q)\triangleq\begin{bmatrix}\pi^\ast_{z_1}(\tilde{\epsilon}(t_q))&\cdots&\pi^\ast_{z_{|\mathcal{Z}|}}(\tilde{\epsilon}(t_q))\end{bmatrix}^T,
\end{align}
where
\begin{align*}
&\pi^\ast_z(\tilde{\epsilon}(t_q))= \frac{\sigma_z(\tilde{\epsilon}(t_q))}{\bm{\sum}\limits_{z'\in \mathcal{Z}}\sigma_{z'}(\tilde{\epsilon}(t_q))}\\
&\sigma_z(\tilde{\epsilon}(t_q))= \bm{\sum}\limits_{T\in G_q(z)}\prod\limits_{(z',z)\in T}P_q^{\tilde{\epsilon}(t_q)}(z',z).
\end{align*}}
\cut{In addition, we assume the observation errors and exploration deviations satisfy the following assumptions.}


\subsection{Main analytical result}\label{sec:main}
The following theorem is the main analytical result of this paper. It shows that the state $z(t)$ converges to the set of stochastically stable state with probability one. 
Moreover, the convergence rate is quantified using the distance between $\pi(t)$ and the limiting distribution $\pi^\ast$; i.e., $D(t)\triangleq ||\pi(t)-\pi^\ast ||$. The formal proof of Theorem~\ref{The4} will be given in Section~\ref{sec:proof}.

\begin{theorem}\label{The4}
Let Assumptions~\ref{asu_error_o} - \ref{asu6} hold, the following properties hold for the RL algorithm:\\
(P1) $\lim\limits_{t\to\infty}Pr\{z(t)\in\Lambda^\ast\}=1$ and $\Lambda^\ast\subseteq {diag} (\mathcal{S}\times\mathcal{S})$;\\
(P2) there exist positive integer $t_{min}$ and positive constant $C$ such that for any $t^\ast>t_{min}$ and $t\ge t^\ast+1$, the following is true:
\begin{align}\label{theorem_in}
%
%
\nonumber&D(t)\le  {C}( ||{\epsilon}(t^\ast)||_\infty+ ||{\epsilon}(t)||_\infty+ e_r(t^\ast)\\
& \exp (-\bm{\sum}\limits_{\tau=t^\ast}^{t-1} \prod\limits_{i=1}^N\epsilon_i(\tau)|\mathcal{A}_i|)+\exp (-\bm{\sum}\limits_{\tau=t^\ast}^{t-1} \prod\limits_{i=1}^N\tilde\epsilon_i(\tau)|\mathcal{A}_i|)).
\end{align}%
%
%
%
\cut{there exists $t^{\ast}_3_\vee >2\mathcal{T}+2$ and $t_e >2\mathcal{T}+2$, for any $t^{\ast}_3>\max\{t^{\ast}_3_{\vee},t_e\}$and $t>t^{\ast}_3$, the convergence rate could be quantified as follows:
\pc{\begin{align*}
&\left| \left| \pi(t)  - \pi^\ast \right|\right|= O\left(\bm{\sum}\limits_{k=t^{\ast}_3}^t |{e}(k)|+\epsilon_\vee(t^{\ast}_3)+\epsilon_\vee(t) \right.\\
&\cut{\quad}\left.+\exp \left(-\bm{\sum}\limits_{\tau=t^{\ast}_3}^{t-\mathcal{T}-1}\epsilon_d(\tau)\epsilon_a(\tau) (nm)\right)\right),
\end{align*}where $t, t^{\ast}_3, k$ and $\tau$ are multiples of $(\mathcal{T}+1)$.}}
%
\end{theorem}

\section {Discussion}\label{Limit_ana}
\subsection{Weakly acyclic games} 
In  this section, we study the special case where the interactions of the players consist of a weakly acyclic game. A game is called to be weakly acyclic if from every action profile, there exists a finite best-response improvement path \cut{(i.e., a finite sequence of players' best responses)}leading from the action profile to a pure Nash equilibrium.\cut{ Then this game is called weakly acyclic game.}\cut{ It is known that any weakly acyclic game has at least one pure Nash equilibrium~\citep{fabrikant2010structure,milchtaich1996congestion,young1993evolution}. } In addition, we further assume that the measurement noises are uniformly upper bounded.
\begin{asu}\label{asu_error_ub}
There is constant $\mathcal{E}\ge 0$, such that $|w_i(t)|\leq \mathcal{E}, \forall i\in\mathcal{V}$ and $t\ge 0$. 
\end{asu}
We will show the convergence of the RL algorithm to the set of equilibria defined as follows  when the game is weakly acyclic.\cut{ when the game is weakly acyclic. First, we give the definition of the special case of pure Nash equilibrium.}

\begin{definition}\label{def1_NE}(Anti-$\mathcal{E}$-Nash equilibrium) An action profile ${s}_\ast\triangleq (a^1_\ast,\cdots,a^i_\ast,\cdots,a^N_\ast)$ is an anti-$\mathcal{E}$-Nash equilibrium if ~$\forall i\in \mathcal{V}, \forall a^i\in\mathcal{A}_i$, $u_i({s}_\ast)\ge u_i(a^i,a^{-i}_\ast)+2\mathcal{E}$.\end{definition}

\begin{remark}\label{rmk1}
Anti-$\mathcal{E}$-Nash equilibrium is stronger than pure Nash equilibrium and $\mathcal{E}$-approximate Nash equilibrium~\cite{nisan2007algorithmic}. When $\mathcal{E}= 0$, anti-$\mathcal{E}$-Nash equilibrium reduces to Nash equilibrium. In the $\mathcal{E}$-approximate Nash equilibrium, the inequality becomes $u_i({s}_\ast)\ge u_i(a^i,a^{-i}_\ast)-2\mathcal{E}$. So any anti-$\mathcal{E}$-Nash equilibrium is also a Nash equilibrium and any Nash equilibrium is an $\mathcal{E}$-approximate Nash equilibrium. The reverses are not true in general. \qed
\end{remark}

Denote the set of anti-$\mathcal{E}$-Nash equilibria of the game $\Gamma$ as $\mathcal{N}_{\mathcal{E}^-}(\Gamma)$ and $diag(\mathcal{N}_{\mathcal{E}^-}(\Gamma)\times\mathcal{N}_{\mathcal{E}^-}(\Gamma))\triangleq\{(s,s)|s\in\mathcal{N}_{\mathcal{E}^-}(\Gamma) \}$. The following corollary implies that the action profiles converge to $\mathcal{N}_{\mathcal{E}^-}(\Gamma)$ with probability one.

\begin{corollary}\label{The2}
If Assumptions~\ref{asu_error_o} - \ref{asu_error_ub}\cut{ and \ref{asu_error}} hold, $\Gamma$ is a weakly acyclic game and $\mathcal{N}_{\mathcal{E}^-}(\Gamma)\neq \emptyset$, then it holds that $\lim \limits_{t\to \infty}Pr\{z(t)\in{diag}(\mathcal{N}_{\mathcal{E}^-}(\Gamma)\times\mathcal{N}_{\mathcal{E}^-}(\Gamma))\}=1$ for the RL algorithm. 
\end{corollary}

From Theorem~\ref{The4}, we have $\lim \limits_{t\to \infty}Pr\{z(t)\in \Lambda^\ast\}=1$ and $\Lambda^\ast\subseteq {diag} (\mathcal{S}\times\mathcal{S})$. Then following the proofs of Lemma 4.2 and Claims 3 -4 in the Proposition 4.3 in~\cite{MZ-SM:SIAM13}, we can get that $\Lambda^\ast\subseteq diag(\mathcal{N}_{\mathcal{E}^-}(\Gamma)\times\mathcal{N}_{\mathcal{E}^-}(\Gamma))$ if $\Gamma$ is weakly acyclic and $\mathcal{N}_{\mathcal{E}^-}(\Gamma)\neq \emptyset$.

\begin{remark}
As shown in \citep{MZ-SM:SIAM13,doi:10.1137/070680199}, when games are weakly acyclic, stochastically stable states are contained in the set of pure Nash equilibrium. Towards our best knowledge, weakly acyclic games are the most general ones which have such property. When a game is not weakly acyclic, stochastically stable states can still be used to characterize where the algorithm converges. So, stochastically stable states are of broader applicability than pure Nash equilibrium. \qed
\end{remark}

\subsection{Estimate of constant $C$ in inequality (\ref{theorem_in}) }
The following corollary estimates constant $C$ in inequality~(\ref{theorem_in}) when the measurement noises are absent. For presentation simplicity, denote $|\gamma|_{min} \triangleq \min\limits_{i\in \mathcal{V}}\gamma_i$,  $C_{min} \triangleq\min\{(\frac{|\gamma|_{min}}{|\mathcal{A}|_{\infty}})^{N|\mathcal{Z}|}, 1\}$, $C_{max} \triangleq\max\{1,||\gamma||_\infty^{N|\mathcal{Z}|}\}$ and $C_\epsilon\triangleq 8 N|\mathcal{Z}|^{|\mathcal{Z}|+4}(N+1)^{|\mathcal{Z}|}2^{N|\mathcal{Z}|} \frac{C_{max} }{C_{min} }$.

\begin{corollary}\label{Cor1}
If Assumptions~\ref{asu_error_o} - \ref{asu6} hold, $w_i(t)=0$ and the exploration rates satisfy that $||\tilde\epsilon(t)||_\infty\le C_{min} /2(N|\mathcal{Z}|^{|\mathcal{Z}|+3}(N+1)^{|\mathcal{Z}|}2^{N|\mathcal{Z}|}\allowbreak C_{max}  )$ for all $t$, then the constant $C$ in inequality~(\ref{theorem_in}) can be estimated as  $C= \max\{|\mathcal{Z}|,4^N, 4C_\epsilon\}$.
\end{corollary}

The proof of Corollary~\ref{Cor1} will be given in Section~\ref{sec:pfCor1}.

\subsection{Explicit convergence rate}
If the exploration rates and exploration deviations are given, we can explicitly quantify how fast the algorithm will reach the set $\Lambda^\ast$.\cut{ In particular, for any small $\delta >0$, we will provide a $t_\delta$ such that $D(t)<\delta, \forall t\ge t_\delta$.} Assume the exploration rate for player $i$ is $\frac{1}{|\mathcal{A}_i|t^{1/N}}$, and the exploration deviations are ${e_i}(t)=\cut{\frac{1}{\tau^{2/N}}} 0,\forall i\in\mathcal{V}$. By Theorem~\ref{The4}, we have:
\begin{align}\label{explicit_in}
\nonumber&D(t)\le {C}\left(2\exp (-\bm{\sum}\limits_{\tau=t^\ast}^{t-1} \frac{1}{\tau})+\frac{1}{ {t^\ast}^{1/N}}+\frac{1}{ t^{1/N}}\right)\\
\nonumber&= {C}\left(2\exp (\bm{\sum}\limits_{\tau=1}^{t^\ast-1} \frac{1}{\tau}-\bm{\sum}\limits_{\tau=1}^{t-1} \frac{1}{\tau})+\frac{1}{ {t^\ast}^{1/N}}+\frac{1}{ t^{1/N}}\right)\\
\nonumber&\le {C}\left(2\exp (1+\int _{1}^{t^{\ast}-1}{\frac{1}{x}}dx-\int _{1}^{t}{\frac{1}{x}}dx)\right.\\
\nonumber&\cut{\quad} \left.+\frac{1}{ {t^\ast}^{1/N}}+\frac{1}{ t^{1/N}}\right)\\
\nonumber&= {C}\left(2\exp \left(1-\int _{t^{\ast}-1}^{t}{\frac{1}{x}}dx\right)+\frac{1}{ {t^\ast}^{1/N}}+\frac{1}{ t^{1/N}}\right)\\
&\le{C}\left(\frac{2e(t^\ast-1)}{t}+\frac{2}{ {t^\ast}^{1/N}}\right).
\end{align}
The second inequality of (\ref{explicit_in}) is a result of inequality (2) of \cite{CHLEBUS2009732}. Given any $\delta>0$, $D(t)\le\delta$ for all \cut{$t^\ast\ge \left(\frac{4C}{\delta}\right)^N$ and} $t\ge \cut{t_\delta=\frac{4Ce(t^\ast-1)}{\delta}=}\frac{e(4C)^{N+1}}{\delta^{N+1}}-\frac{4Ce}{\delta}$. Roughly speaking, it takes $O(\frac{1}{\delta^{N+1}})$ iterations to reach error $\delta$.

\subsection{Optimal exploration rates}\label{sec:optexp}
An interesting question is how to choose the exploration rates to minimize the upper bound in inequality (\ref{theorem_in}). This is an infinite-dimension and non-convex optimization problem and hard to solve in general. For analytical tractability, we restrict the exploration rates to be $\mathbf{p}$-series which have been widely used in stochastic approximation and convex optimization~\citep{bertsekas2015convex,hasʹminskii1972stochastic,citeulike:2621242}. In particular, let ${e_i}(t)=\cut{\frac{1}{\tau^{2/N}}} 0$ and  $\epsilon_i(t)=\frac{1}{|\mathcal{A}_i|t^{\mathbf{p}/N}}, \mathbf{p}\in(0,1], \forall i\in\mathcal{V}$. This choice satisfies Assumption \ref{asu6}. We aim to choose $\mathbf{p}\in(0,1]$ to minimize the upper bound of $D(t)$. With such restriction, inequality (\ref{theorem_in}) becomes:
\begin{align}\label{explicit_in_o}
\nonumber&D(t)\le {C}(2\exp (-\bm{\sum}\limits_{\tau=t^\ast}^{t-1} \frac{1}{\tau^\mathbf{p}})+\frac{1}{ {t^\ast}^{\mathbf{p}/N}}+\frac{1}{ t^{\mathbf{p}/N}})\\
&\le {C}(2\exp (1-\int _{t^{\ast}-1}^{t}{\frac{1}{x^\mathbf{p}}}dx)+\frac{1}{ {t^\ast}^{\mathbf{p}/N}}+\frac{1}{ t^{\mathbf{p}/N}}).
\end{align}
The second inequality of (\ref{explicit_in_o}) follows the same steps of (\ref{explicit_in}) by replacing $\frac{1}{\tau}$ with $\frac{1}{\tau^\mathbf{p}}$.
When $\mathbf{p}\in(0,1)$, inequality (\ref{explicit_in_o}) becomes:
\begin{align*}
&D(t) \\
&\le{C}(2\exp (1+\frac{(t^\ast-1)^{1-\mathbf{p}}-t^{1-\mathbf{p}}}{1-\mathbf{p}})+\frac{1}{ {t^\ast}^{\mathbf{p}/N}}+\frac{1}{ t^{\mathbf{p}/N}})\\
&= {C}\left(2\exp (1+\frac{(t^\ast-1)^{1-\mathbf{p}}}{1-\mathbf{p}})\exp(\frac{-t^{1-\mathbf{p}}}{1-\mathbf{p}})\right.\\
&\cut{\quad}\left.+\frac{1}{ {t^\ast}^{\mathbf{p}/N}}+\frac{1}{ t^{\mathbf{p}/N}}\right).
\end{align*}
Since $\lim\limits_{t\to \infty}\frac{\frac{t^{1-\mathbf{p}}}{1-\mathbf{p}}}{\frac{\mathbf{p}}{N}\ln t}\allowbreak=\infty$, we have $\lim\limits_{t\to \infty}\frac{\exp(\frac{-t^{1-\mathbf{p}}}{1-\mathbf{p}})}{\exp(-\frac{\mathbf{p}}{N}\ln t)}\allowbreak=0\cut{<\infty}$. So the term $\frac{1}{ t^{\mathbf{p}/N}}$ dominates the term ${C}(2\exp (1+\frac{(t^\ast-1)^{1-\mathbf{p}}}{1-\mathbf{p}})\exp(\frac{-t^{1-\mathbf{p}}}{1-\mathbf{p}})$ as $t$ increases. 
When $\mathbf{p}=1$, inequality (\ref{explicit_in_o}) becomes: 
\begin{align*}
&D(t)\le  {C}(2\exp (1+\ln(t^\ast-1)-\ln t)+\frac{1}{ {t^\ast}^{1/N}}+\frac{1}{ t^{1/N}})\\
&= {C}(2\exp (1+\ln(t^\ast-1)) \frac{1}{t}+\frac{1}{ {t^\ast}^{1/N}}+\frac{1}{ t^{1/N}}).
\end{align*}
Analogously, we have $\lim\limits_{t\to \infty}(\frac{1}{t})/(\frac{1}{ t^{1/N}})=0\cut{<\infty}$.  So the term $\frac{1}{ t^{1/N}}$ dominates the term $\frac{1}{t}$ as $t$ increases. In both cases, $\frac{1}{ t^{\mathbf{p}/N}}$ dominates the upper bound in (\ref{explicit_in_o}).
\cut{Notice that the $\frac{1}{ t^{\mathbf{p}/N}}$ is minimized by choosing $\mathbf{p}=1$.}When $\mathbf{p}=1$, $\frac{1}{ t^{\mathbf{p}/N}}$ decreases fastest among $\mathbf{p}\in(0,1]$. Therefore, $\epsilon_i(t)=\frac{1}{|\mathcal{A}_i|t^{1/N}}$ is optimal among $\mathbf{p}$-series.

\cut{Then The problem can be formulated as follows:
\begin{align}\label{Problem_A}
\nonumber\min\limits_{a\in \mathbb{R}} &J_j({\epsilon}(\tau))\triangleq\exp (-\bm{\sum}\limits_{\tau\in\{t^{\ast}_3+\ell^\ast\}} \tau^{aN}\prod\limits_{i=1}^N\lambda_i|\mathcal{A}_i|)\\
\nonumber& + ||{\epsilon}(t^{\ast}_3)||_\infty+ |{\epsilon}(j)|_\infty+ \frac{{t^{\ast}_3}^{-2N}} {\prod\limits_{i=1}^N(\lambda_i {t^{\ast}_3}^a+{t^{\ast}_3}^{-2})} \\
\nonumber&+\exp (-\bm{\sum}\limits_{\tau\in\{t^{\ast}_3+\ell^\ast\}} \prod\limits_{i=1}^N(\lambda_i \tau^a+\tau^{-2})|\mathcal{A}_i|)\\
\nonumber\text{subject to } &\ell^\ast = 0, ..., j-1\tag{PA}\\
\nonumber& \lim\limits_{j\to\infty}||{\epsilon}(t^{\ast}_3)||_\infty=0\\
\nonumber& \lim\limits_{j\to\infty}\frac{{j}^{-2N}} {\prod\limits_{i=1}^N(\lambda_i {j}^a+{j}^{-2})}=0\\
\nonumber& \bm{\sum}\limits_{\tau\in\{t^{\ast}_3+\ell^\ast\}}\prod\limits_{i=1}^N(\lambda_i \tau^a+\tau^{-2})=\infty\\ 
\nonumber& \bm{\sum}\limits_{\tau\in\{t^{\ast}_3+\ell^\ast\}}\tau^{aN}\prod\limits_{i=1}^N\lambda_i=\infty.
\end{align}}

\cut{
\cut{ From Theorem~\ref{The4}, we can see the convergence speed is affected by both the attacker and defender. }\cut{ Even if the defender cannot totally control the convergence speed, the convergence rate can still provide some meaningful prediction in real-world implementation. And if given some exploration rates of the attacker and exploration deviations, the defender may choose its own exploration rates to make the convergence faster. }If the exploration rates, utility delays and exploration deviations are given, we can explicitly quantify how fast the algorithm will reach the set of action profiles with maximum stochastic potential. For example, if the exploration rates $\epsilon_d(t), \epsilon_a(t)$ are both $\frac{1}{\sqrt{tnm}}$, the exploration deviations $e^d_c(t), e^a_c(t)$ are both $0$ and the utility delays are $0$, then \begin{align*}
D(j)\le \exp\left(-\int _{t^{\ast}_3}^{j}{\frac {1}{x}}dx\right)+\frac{C}{\sqrt{j}} \le \frac{t^{\ast}_3}{j}+\frac{C}{\sqrt{j}}.
\end{align*}\cut{That is $D(j)$ goes to $0$ as $j$ goes to infinity and }The convergence rate is at the order of $\frac{1}{\sqrt{j}}$.

}


\subsection{Memory and communication}
The RL algorithm only requires each player to remember its own utility values and actions in recent history. So the memory cost is low. In addition, communications are case dependent. In \cite{MZ-SM:SIAM13}, the utility function of each robot only depends on the actions of its own and nearby robots. The communication range of each robot is twice of its sensing range. So the communication graphs are time-varying and usually sparse. In Section~\ref{sec:case1}, each customer can communicate with the system operator. So the communication graph is a fixed star graph.


\cut{
{\bf Multiple-player game.} The adaptive algorithm can be extended to the multiple-player game. If each player uses the similar exploration- exploitation strategy, the action profile can converge to the set of profiles with maximum stochastic potential. And if the game is a weakly acyclic game, the action profile will converge to the set of Nash equilibria. The multiple-player game can be applied to coverage optimization problem and related discussed in~\cite{MZ-SM:SIAM13}.}

\cut{We still consider the Markov chain $\{P_t\}$ introduced in Section~\ref{Notation1}. Bbecause we consider the errors, the transition matrices at different $t$ might be different from $P^{{\epsilon}(t)}$. Here we define the transition matrices as $\tilde P^{{\epsilon}(t)}$. In Theorem~\ref{The1}, we actually prove that not matter what the initial distribution $\pi(0)$ is, $\pi(t)=\pi(0)P^{{\epsilon}(1)}\cdots P^{{\epsilon}(t)}$ will converge to the limiting distribution $\pi^\ast$ shown in Lemma~\ref{lemma:star2}. Now if for any initial distribution $\pi(0)$, we can prove that $\pi(t)=\pi(0)\tilde P^{{\epsilon}(1)}\cdots \tilde P^{{\epsilon}(t)}$ will converge to $\pi^\ast$ as well, we can get the same asymptotic convergence result as Theorem~\ref{The1}. To prove this result, we need the following two lemmas.}

\cut{
\begin{lemma}\label{lemma:star3}
If the sequence of transition matrices $\tilde P^{{\epsilon}(t)}$ satisfies that $\bm{\sum}\limits_{t=1}^\infty ||\tilde P^{{\epsilon}(t)}-P^{{\epsilon}(t)}||< \infty$, where the norm of matrix is $||P||\triangleq \max\limits_{i}\bm{\sum}\limits_j|p_{ij}|$. Then for any initial distribution $\pi(0)$, $\lim\limits_{t\to \infty}\pi(0)\tilde P^{{\epsilon}(1)}\cdots \tilde P^{{\epsilon}(t)}=\pi^\ast$.
\end{lemma}

The proof of this lemma is given in Lemma 3.1 of \cite{gutjahr1996simulated}.
}
\cut{
\pc{
By Lemma~\ref{lemma:star}, the support of $\pi^\ast$ is contained in $\Lambda^\ast$. Therefore we have $\lim \limits_{t\to \infty}Pr\{(s(t-1),s(t))\in \Lambda^\ast\}=1$. Now we want to prove that $\Lambda^\ast\subseteq diag(\mathcal{S}_{BR}(a^\ast)\times \mathcal{S}_{BR}(a^\ast))$.

The proof is divided into 3 claims.

{\bf Claim 1.} There exist a \zs{$t_F\le\hat{t}<\infty$}. For all $t\ge \hat{t}$, $\Lambda({\epsilon}(t))\subseteq {diag} (\mathcal{S}\times\mathcal{S})$.
\begin{pf}
\zs{For all $t\ge t_F$,} assume $({s}_0,{s}_1)\in \Lambda({\epsilon}(t))$ but $({s}_0,{s}_1)\not\in {diag} (\mathcal{S}\times\mathcal{S})$, i.e., ${s}_0\neq {s}_1$, where ${s}_0=(d^0,a^0),{s}_1=(d^1,a^1)$. Since $({s}_0,{s}_1)\in \Lambda({\epsilon}(t))$, then there is a tree $T_{max}({\epsilon}(t))$ rooted at $({s}_0,{s}_1)$ such that it has largest total probability. Here we consider two cases: {\bf case 1.} $\mu_{d^0}\ge\mu_{d^1}$; {\bf case 2.} $\mu_{d^1}>\mu_{d^0}$.

For case 1, we construct a tree $T'$ by adding some feasible edges from $({s}_0,{s}_1)$ to $(\hat{s},\hat{s})$,where $\hat{s}=(d^0,a^1)$, and deleting the edge of leaving $(\hat{s},\hat{s})$. For $t\ge t_F$, $\mu_{d^0}\ge\mu_{d^1}$ implies $\bar\mu_{d^0}(t)\ge\bar\mu_{d^1}(t)$. \zs{The transition $({s}_0,{s}_1)\to({s}_1,\hat{s})$ is achieved when both of the defender and attacker perform exploitation. Because $\bar\mu_{d^0}\ge \bar\mu_{d^1}$, by Line 20 in Algorithm~\ref{algo:adaptive}, the defender chooses $d^0$. And the transition $({s}_0,{s}_1)\to({s}_1,\hat{s})$ is also achieved when both of the defender and attacker perform exploitation. By equalities (2)-(5), A possible path from $({s}_0,{s}_1)$ to $(\hat{s},\hat{s})$ is}:
\begin{align*}
&({s}_0,{s}_1)\xrightarrow {(1-\epsilon_d(t))(1-\epsilon_a(t))} ({s}_1,\hat{s})\\
&({s}_1,\hat{s}) \xrightarrow {(1-\epsilon_d(t))(1-\epsilon_a(t))} (\hat{s},\hat{s}),
\end{align*}

Then the product of the probabilities of the added edges from $((d^0,a^0),(d^1,a^1))$ to $((d^0,a^1),(d^0,a^1))$ is $((1-\epsilon_d(t))(1-\epsilon_a(t)))^2$.

Let us say the next step of $(\hat{s},\hat{s})$ is $(\hat{s},{s}_2)$, where ${s}_2=(d^2,a^2)$ and $d^2 ,a^2$ satisfy that $d^2\neq d^0,a^2=a^1$. Then the transition probability of the leaving edge can be represented as follows.
\begin{align*}
(\hat{s},\hat{s})\xrightarrow {(1-\epsilon_a(t))\frac{\epsilon_d(t)}{n-2•} }(\hat{s},{s}_2).
\end{align*}
Since $(1-\epsilon_a(t))\frac{\epsilon_d(t)}{n-2•}<\epsilon_d(t)$, then the transition probability of the deleted edge is less than $\epsilon_d(t)$.

For case 2, we construct a tree $T'$ by adding some feasible edges from $({s}_0,{s}_1)$ to $({s}_1,{s}_1)$, and deleting the edge of leaving $({s}_1,{s}_1)$. Similar to case 1, a possible transition from $({s}_0,{s}_1)$ to $({s}_1,{s}_1)$ is:
\begin{align*}
&({s}_0,{s}_1)\xrightarrow {(1-\epsilon_d(t))(1-\epsilon_a(t))} ({s}_1,{s}_1).
\end{align*}
Similar to case 1, the transition probability of the deleted edge is less than $\epsilon_d(t)$.

By (1) and (2) in Assumption~\ref{asu6}, $\epsilon_a(t), \epsilon_d(t)$ are strictly decreasing to $0$. Hence, there exists a \zs{$t_F\le\hat{t}<\infty$} such that $\forall t\ge \hat{t},\epsilon_d(t)< ((1-\epsilon_d(t))(1-\epsilon_a(t)))^2\le (1-\epsilon_d(t))(1-\epsilon_a(t))$.
So for both cases, the total probability of $T'$ is larger than that of $T_{max}({\epsilon}(t))$. We reach a contradiction.
\qed\end{pf}

{\bf Claim 2.} There exists a \zs{$t_F\le t_0<\infty$}, for all $t\ge \max\{\hat{t},t_0\}$, $\Lambda({\epsilon}(t))\subseteq {diag} (\mathcal{S}_{BR}(a(t))\times\mathcal{S}_{BR}(a(t)))$.
\begin{pf}
By Claim 1, for any $t\ge \hat{t}$, we have $\Lambda({\epsilon}(t))\subseteq {diag} (\mathcal{S}\times\mathcal{S})$. Consider $t\ge \hat{t}$, assume $(s,s)\in \Lambda({\epsilon}(t))$, where $s=(d,a)$. But $(s,s)\not\in {diag} (\mathcal{S}_{BR}(a)\times\mathcal{S}_{BR}(a))$. And since $(s,s)\in \Lambda({\epsilon}(t))$, then there is a tree $T_{max}({\epsilon}(t))$ rooted at $(s,s)$ such that it has largest total probability. Construct a tree $T'$ by adding the feasible edge from $(s,s)$ to $({s}_\ast,{s}_\ast)$ and deleting the edge of leaving $({s}_\ast,{s}_\ast)$, where ${s}_\ast\in S_{BR}(a)$, ${s}_\ast=(d^\ast,a)$ and $d^\ast \neq d$. According to Line 19 - 27 in Algorithm~\ref{algo:adaptive}, the path from $(s,s)$ to $({s}_\ast,{s}_\ast)$ is the case as follows.

\begin{align*}
&(s,s) \xrightarrow {(1-\epsilon_a(t))\frac{\epsilon_d(t)}{n-2•}} (s,{s}_\ast)\\
&(s,{s}_\ast) \xrightarrow {(1-\epsilon_d(t))(1-\epsilon_a(t))} ({s}_\ast,{s}_\ast).
\end{align*}
Then the probability of the added edge is
\begin{align*}
&\left((1-\epsilon_a(t))\frac{\epsilon_d(t)}{n-2•}\right)(1-\epsilon_d(t))(1-\epsilon_a(t)).
\end{align*}
Consider the edge of leaving $({s}_\ast,{s}_\ast)$. Let us say the next step of $({s}_\ast,{s}_\ast)$ is $({s}_\ast,{s}_2)$, where ${s}_2=(d^2,a^2)$ and $d^2 ,a^2$ satisfy that $d^2\neq d^\ast$ and $a^2\neq a$. Then the probability of the deleted edge is $\frac{\epsilon_d(t)\epsilon_a(t) }{(n-2)(m-2)}$.

By (1) and (2) in Assumption~\ref{asu6}, $\epsilon_a(t), \epsilon_d(t)$ are strictly decreasing to $0$. Hence, there exists a \zs{$t_F\le t_0<\infty$} such that $\forall t\ge\max\{\hat{t},t_0\}$,
\begin{align*}
&\frac{\epsilon_d(t)\epsilon_a(t) }{(n-2)(m-2)}<\\
&\left((1-\epsilon_a(t))\frac{\epsilon_d(t)}{n-2•}\right)(1-\epsilon_d(t))(1-\epsilon_a(t)).
\end{align*}
Therefore, the total probability of $T'$ is larger than that of $T_{max}({\epsilon}(t))$. We reach a contradiction.
\qed\end{pf}

{\bf Claim 3.} $\Lambda^\ast\subseteq {diag} (\mathcal{S}_{BR}(a^\ast)\times\mathcal{S}_{BR}(a^\ast))$.
\begin{pf}
Claim 2 holds for any ${\epsilon}(t)$ if $t\ge \max\{\hat{t},t_0\}$, and by 2) in Assumption~\ref{asu6}, $\lim \limits_{t\to \infty}{\epsilon}(t)=\epsilon^\ast$, and $\lim \limits_{t\to \infty}a(t)=a^\ast$. Then we reach the desired result.
\qed\end{pf}
}
}

\section{Proofs}\label{sec:proof}
We will prove Theorem~\ref{The4} and Corollary~\ref{Cor1} in this section. 
\subsection{Analysis of the H-RL algorithm}
For the sake of analysis, we introduce the H-RL algorithm, which has time-homogeneous exploration rates; i.e., ${\tilde{\epsilon}(t)}={\tilde{\epsilon}}\in(0,1], \forall t\ge 0$ in the RL algorithm. Then $\{z(t)\}$ in the H-RL algorithm forms a time-homogeneous Markov chain $\mathcal{M}^{\tilde{\epsilon}}$ withbe the transition matrix $P^{\tilde{\epsilon}}$. The analysis of the H-RL algorithm provides preliminary results for that of the RL algorithm.
The following lemma studies the properties of the feasible transitions in the Markov chain $\mathcal{M}^{\tilde{\epsilon}}$.

\begin{lemma}\label{lemma:nonzero2}
Given any ${\tilde{\epsilon}}\in(0,1]$, if Assumption \ref{asu_error_o} holds, then each nonzero entry in transition matrix $ P^{\tilde{\epsilon}}$ is a polynomial of the variables $\{\tilde\epsilon_i, 1-\tilde\epsilon_i\}_{i\in \mathcal{V}}$. In addition, the coefficients of the polynomials are independent of $\tilde\epsilon$.
\end{lemma}

\begin{pf}
Consider any two states $x,y \in \mathcal{Z}$ that the transition from $x$ to $y$ is feasible within one step. In particular, $x=({s}(0),{s}(1))$ and $y=({s}(1),{s}(2))$, where ${s}(t)=(a^1(t),\cdots,a^N(t))$ for $t\in\{0,1,2\}$. And the transition probability is $P^{\tilde{\epsilon}}(x,y)=\prod\limits_{ i\in \mathcal{V}}Pr\{a^i(2)|a^i(0),a^i(1)\}$.

If player $i$ performs exploitation with $a^{i}(0)\neq a^{i}(1)$ and $u_{i}({s}(0))\ge u_{i}({s}(1))$ but $a^{i}(2)=a^{i}(1)$, we call this event a mis-exploitation between ${s}(0)$ and ${s}(1)$. By Assumption~\ref{asu_error_o}, the mis-exploitation between ${s}(0)$ and ${s}(1)$ of player ${i}$ happens with probability $\delta_{i}({s}(0),{s}(1))\triangleq Pr\{w_i(0)-w_i(1)<u_{i}({s}(1))-u_{i}({s}(0))\}$ if $a^{i}(0)\neq a^{i}(1)$ and $\delta_{i}({s}(0),{s}(1))\triangleq 0$ if $a^{i}(0)= a^{i}(1)$.
%
The probability $\delta_{i}({s}(0),{s}(1))$ is independent of $\tilde\epsilon$ and only determined by the two utility values $u_{i}({s}(0))$ and $u_{i}({s}(1))$ because $\{w_i(t)\}_{t\ge 0}$ are I.I.D. Given states $x$ and $y$, the set of players can be partitioned into three sets: $ \mathcal{V}_{er}(x,y)\triangleq \{i\in\mathcal{V}|a^i(2)\not\in \{a^i(0),a^i(1)\}\}$, $ \mathcal{V}_{ex}(x,y)\triangleq \{j\in\mathcal{V}|a^j(2)=a^j(\arg\max\limits_{t\in\{0,1\}}\{u_{j}({s}(t))\})\}$ and $ \mathcal{V}_{m}(x,y)\triangleq \{k\in\mathcal{V}|\delta_{k}({s}(0),{s}(1))\neq 0\}\cap\{k\in\mathcal{V}|a^k(2)=a^k(\arg\min\limits_{t\in\{0,1\}}\{u_{k}({s}(t))\})\}$. For any $ i\in \mathcal{V}_{er}(x,y), Pr\{a^i(2)|a^i(0),a^i(1)\}=\frac{\tilde\epsilon_i}{|\mathcal{A}_i|}$. And for any $j\in \mathcal{V}_{ex}(x,y)$, $a^j(2)$ can be achieved by exploitation without mis-exploitation or exploration, then $Pr\{a^j(2)|a^j(0),a^j(1)\}=(1-\delta_j({s}(0),{s}(1)))(1-\tilde\epsilon_j) + \frac{\tilde\epsilon_j}{|\mathcal{A}_j|}$. Similarly, for any $k\in \mathcal{V}_{m}(x,y)$, $a^k(2)$ can be achieved by mis-exploitation or exploration, then $Pr\{a^k(2)|a^k(0),a^k(1)\}\allowbreak=\delta_k({s}(0),{s}(1))(1-\tilde\epsilon_k) + \frac{\tilde\epsilon_k}{|\mathcal{A}_k|}$. 
Then the transition probability can be written as: 
\begin{align}\label{Lem1e1}
\nonumber&P^{\tilde{\epsilon}}(x,y)=\prod\limits_{ i\in \mathcal{V}_{er}(x,y)}\frac{\tilde\epsilon_i}{|\mathcal{A}_i|}\\
\nonumber&\times\prod\limits_{j\in \mathcal{V}_{ex}(x,y)}\left((1-\delta_j({s}(0),{s}(1)))(1-\tilde\epsilon_j) + \frac{\tilde\epsilon_j}{|\mathcal{A}_j|}\right)\\
&\times\prod\limits_{k\in \mathcal{V}_{m}(x,y)}\left(\delta_k({s}(0),{s}(1))(1-\tilde\epsilon_k) + \frac{\tilde\epsilon_k}{|\mathcal{A}_k|}\right).
\end{align} It is clear that $P^{\tilde{\epsilon}}(x,y)$ is a polynomial of $\{\tilde\epsilon_i, 1-\tilde\epsilon_i\}_{i\in \mathcal{V}}$ with coefficients independent of $\tilde\epsilon$.\qed\end{pf}

Given any ${\tilde{\epsilon}}\in(0,1]$, define stochastic vector $\pi^\ast({\tilde{\epsilon}})$ as the stationary distribution of the Markov chain $\mathcal{M}^{\tilde{\epsilon}}$; i.e., ${ \pi^\ast({\tilde{\epsilon}}) }^T P^{\tilde{\epsilon}}={ \pi^\ast({\tilde{\epsilon}}) }^T$. In the H-RL algorithm, when a player performs exploration, it can choose any element in its action set. One can see that, for any pair of states $x,y\in \mathcal{Z}$, $y$ can be reached from $x$ within finite steps, and Markov chain $\mathcal{M}^{\tilde{\epsilon}}$ is ergodic. By Lemma 3.1 in Chapter 6 of \cite{freidlin2012random}, $\pi^\ast({\tilde{\epsilon}})$ can be written as follows:
\begin{align}
\pi^\ast({\tilde{\epsilon}})\triangleq\begin{bmatrix}\pi^\ast_{z_1}(\tilde{\epsilon})&\cdots&\pi^\ast_{z_{|\mathcal{Z}|}}(\tilde{\epsilon})\end{bmatrix}^T,\label{S4e1}
\end{align}
where $\pi^\ast_z(\tilde{\epsilon})= \frac{\sigma_z(\tilde{\epsilon})}{\bm{\sum}\limits_{z'\in \mathcal{Z}}\sigma_{z'}(\tilde{\epsilon})}, \sigma_z(\tilde{\epsilon})= \bm{\sum}\limits_{T\in G_{\tilde{\epsilon}}(z)}\prod\limits_{(z',z)\in E(T)}\allowbreak P^{\tilde{\epsilon}}(z',z)$ and $E(T)$ is the edge set of tree $T$.

\subsection{Stationary distributions without exploration deviations}

Now let us consider an auxiliary scenario that exploration deviations ${e_i}(t)=0$ for all $t$ and for all $i\in\mathcal{V}$. \cut{Given any $t$, we got a fixed ${{\epsilon}(t)}\in(0,1]$. Denote by $P^{{\epsilon}(t)}$ the transition matrix of the time-inhomogeneous Markov chain $\mathcal{M}$ induced by the RL algorithm at iteration $t$. }Then stochastic vector $\hat\pi^\ast({{\epsilon}(t)})$ such that ${ \hat\pi^\ast({{\epsilon}(t)}) }^T P^{{\epsilon}(t)}={ \hat\pi^\ast({{\epsilon}(t)}) }^T$ has the same form of (\ref{S4e1}) with exploration deviations being $0$. For notational simplicity, we refer to $\hat\pi^\ast({\epsilon}(t))$ as $\hat\pi^\ast(t)$. The following lemma shows that $\{\hat\pi^\ast(t)\}_{t\ge 0}$ converges to a limiting distribution with a certain rate, and the support of the limiting distribution is $\Lambda^\ast$.



\begin{lemma}\label{lemma:star2}
If Assumptions~\ref{asu_error_o} - \ref{asu6} hold and ${e_i}(t)=0$ for all $t$ and all $i$, then the sequence $\{\hat\pi^\ast(t)\}_{t\ge 0}$ converges to the limiting distribution $\pi^\ast$ whose support is $\Lambda^\ast\subseteq {diag} (\mathcal{S}\times\mathcal{S})$. Moreover, the convergence rate could be quantified as: $\left| \left|\hat\pi^\ast(t) - \pi^\ast \right|\right|\le C_\epsilon ||{\epsilon}(t)||_\infty$\cut{ for all $t\ge t_e$} for some constant $C_\epsilon>0$.
\end{lemma}
\begin{pf}
The proof is divided into two claims.
\end{pf}

{\bf Claim 1. The limiting distribution $\pi^\ast\triangleq \lim\limits_{t\to\infty}\hat\pi^\ast(t)$ exists, and its support is $\Lambda^\ast\subseteq {diag} (\mathcal{S}\times\mathcal{S})$}.
\begin{pf}
By Lemma~\ref{lemma:nonzero2}, for any ${\epsilon}(t)\in(0,1]$, the non-zero entries of $ P^{{\epsilon}(t)}$ are polynomials of $\{\epsilon_i(t), 1-\epsilon_i(t)\}$ since ${e_i}(t)=0$ for all $i\in \mathcal{V}$ with time-homogeneous coefficients. Then $\sigma_z({\epsilon}(t))$ and $\bm{\sum}\limits_{z'\in \mathcal{Z}}\sigma_{z'}({\epsilon}(t))$ are polynomials of $\{\epsilon_i(t), 1-\epsilon_i(t)\}$ with time-homogeneous coefficients. Recall that $\epsilon_i(t)=\gamma_i{\epsilon}^c(t)$.
For particular state $z\in\mathcal{Z}$, $\sigma_z({\epsilon}(t))$ and $\bm{\sum}\limits_{z'\in \mathcal{Z}}\sigma_{z'}({\epsilon}(t))$ are polynomials of ${\epsilon}^c(t)$, and $\hat\pi^\ast_z({\epsilon}(t))$ is a ratio of two polynomials of ${\epsilon}^c(t)$:
\begin{align}\label{L2c1eq1}
\hat\pi^\ast_z({\epsilon}(t))=\frac{\alpha_z({\epsilon}^c(t))}{\beta({\epsilon}^c(t))•}.
\end{align}
In particular,
\begin{align*}
&\alpha_z({\epsilon}^c(t))=b^z_k{\epsilon}^c(t)^k+b^z_{k+1}{\epsilon}^c(t)^{k+1}+\cdots+b^z_h{\epsilon}^c(t)^h\\
&\beta({\epsilon}^c(t))= b_k{\epsilon}^c(t)^k+b_{k+1}{\epsilon}^c(t)^{k+1}+\cdots+b_h{\epsilon}^c(t)^h,
\end{align*}
where $k\ge 0$. Without loss of generality, we assume that of $b_k$ is non-zero. 
When ${\epsilon}^c(t)$ is sufficiently small, $b^z_k{\epsilon}^c(t)^k$ and $b_k{\epsilon}^c(t)^k$ dominate $\frac{\alpha_z({\epsilon}^c(t))}{\beta({\epsilon}^c(t))•}$. Then the limit of the $\hat\pi^\ast_z({\epsilon}(t))$ can be represented as $\pi^\ast_z=\lim\limits_{t\to \infty}\hat\pi^\ast_z({\epsilon}(t))=\frac{b^z_k}{b_k•}$. Note that $b^z_k$ and $b_k$ are also time-homogeneous. 

By the definition of $\Lambda^\ast$, if $b^z_k$ is non-zero, then $z\in \Lambda^\ast$\cut{ and $\pi^\ast_z$ is non-zero}; And if $b^z_k=0$, then $z\not\in \Lambda^\ast$\cut{, which implies $\pi^\ast_z=0$}. We know $\pi^\ast_z=\frac{b^z_k}{b_k•}$, therefore the support of $\pi^\ast$ is contained in $\Lambda^\ast$. 

For any ${\epsilon}(t)\in(0,1]$ Assume a state $z=({s}(0),{s}(1))\in \Lambda({\epsilon}(t))$ but $z\not\in {diag} (\mathcal{S}\times\mathcal{S})$, i.e., ${s}(0)\neq {s}(1)$, where ${s}(0)=(a^1(0),\cdots,a^N(0)),\allowbreak {s}(1)=(a^1(1),\cdots,a^N(1))$. Since $z\in \Lambda({\epsilon}(t))$, then there is a tree $T_{max}({\epsilon}(t))$ rooted at $z$ such that it has largest total probability. We construct a tree $T'$ by adding the following path from $z$ to $z'=({s}(2),{s}(2))$, where $\hat z=({s}(1),{s}(2))$, where ${s}(2)=(a^1(2),\cdots,a^N(2))$ and $a^i(2)=a^i(\arg\max\limits_{\tau\in\{0,1\}}\{u_{i}({s}(\tau))\}), \forall i$:
\begin{align*}
 z \xrightarrow {P^{{\epsilon}(t)}(z,\hat z)}\hat z\xrightarrow {P^{{\epsilon}(t)}(\hat z,z')}z'.
\end{align*}
By Lemma~\ref{lemma:nonzero2}, $\mathcal{V}_{ex}(z,\hat z)=\mathcal{V}_{ex}(\hat z,z')=\mathcal{V}, \mathcal{V}_{er}(z,\hat z)=\mathcal{V}_{er}(\hat z,z')=\emptyset$ and $\mathcal{V}_{m}(z,\hat z)=\mathcal{V}_{m}(\hat z,z')=\emptyset$. So $P^{{\epsilon}(t)}(z,\hat z)=\prod\limits_{i\in \mathcal{V}}\left((1-\delta_i({s}(0),{s}(1)))(1-\epsilon_i(t)) + \frac{\epsilon_i(t)}{|\mathcal{A}_i|}\right)$ and \\$P^{{\epsilon}(t)}(\hat z,z')=\prod\limits_{i\in \mathcal{V}}\left((1-\delta_i({s}(1),{s}(2)))(1-\epsilon_i(t)) + \frac{\epsilon_i(t)}{|\mathcal{A}_i|}\right)$.

Let us consider the edge leaving $z'$:
\begin{align*}
z'\xrightarrow {P^{{\epsilon}(t)}(z',z'')}\hat z''=({s}(2),{s}(3)),
\end{align*}
where ${s}(3)=(a^1(3),\cdots,a^N(3))$ with at least one player $i$ such that $a^i(3)\neq a^i(2)$. That is $|\mathcal{V}_{er}(z',z'')|\ge 1$. Then the transition probability of the leaving edge satisfies:\begin{align*}
&P^{{\epsilon}(t)}(z',z'')= \prod\limits_{ j\in \mathcal{V}_{er}(z',z'')}\frac{\epsilon_j(t)}{|\mathcal{A}_j|}\\
&\prod\limits_{i\in \mathcal{V}_{ex}(z',z'')}\left((1-\delta_i({s}(2),{s}(2)))(1-\epsilon_i(t)) + \frac{\epsilon_i(t)}{|\mathcal{A}_i|}\right).
\end{align*}
By Assumption~\ref{asu_error_o}, $\delta_i({s},{s}')$ is non-negative constant and independent of ${\epsilon}(t)$ for all $i\in\mathcal{V}$ and $s,s'\in\mathcal{S}$. Then $P^{{\epsilon}(t)}(z,\hat z),P^{{\epsilon}(t)}(\hat z,z')$ and $P^{{\epsilon}(t)}(z', z'')$ are dominated by their lowest degree terms when ${\epsilon}(t)$ is sufficiently small. In particular, the lowest degree terms of $P^{{\epsilon}(t)}(z,\hat z)$ and $P^{{\epsilon}(t)}(\hat z,z')$ are constant terms while the lowest degree term $P^{{\epsilon}(t)}(z', z'')$ is at least first-degree term.
Then there exists some $\hat\epsilon_M \in (0,1]$ such that $P^{{\epsilon}(t)}(z,\hat z)P^{{\epsilon}(t)}(\hat z,z') > P^{{\epsilon}(t)}(z', z''), \forall {\epsilon}(t)\in(0,{\epsilon}_M)$. That is, the total probability of $T'$ is larger than that of $T_{max}({\epsilon}(t))$. We reach a contradiction. Therefore, $\Lambda^\ast\subseteq {diag} (\mathcal{S}\times\mathcal{S})$ because $\Lambda({\epsilon}(t))\subseteq {diag} (\mathcal{S}\times\mathcal{S})$ holds for any sufficiently small ${\epsilon}(t)$. \qed
\end{pf}
{\bf Claim 2. $\left| \left| \hat\pi^\ast(t) - \pi^\ast \right|\right|\le  C_\epsilon ||{\epsilon}(t)||_\infty$ for some constant  $C_\epsilon>0$.}\cut{ for $t\ge t_e$}
\begin{pf}
From Claim 1, we have  $\lim\limits_{t\to\infty}\hat\pi^\ast(t)=\pi^\ast$ and the support of $\pi^\ast$ is $\Lambda^\ast$. Now we consider the convergence rate. 
\begin{align}\label{L2c2eq1}
\nonumber &\left| \left|\hat\pi^\ast(t) - \pi^\ast \right|\right|=\bm{\sum}\limits_{z\in\mathcal{Z}}|\hat\pi^\ast_z({\epsilon}(t))-\pi^\ast_z|\\
\nonumber &=\bm{\sum}\limits_{z\in\Lambda^\ast}|\hat\pi^\ast_z({\epsilon}(t))-\pi^\ast_z|+\bm{\sum}\limits_{z\not\in\Lambda^\ast}|\hat\pi^\ast_z({\epsilon}(t))-\pi^\ast_z|\\
\nonumber&=\bm{\sum}\limits_{z\in\Lambda^\ast}|\hat\pi^\ast_z({\epsilon}(t))-\pi^\ast_z|+\bm{\sum}\limits_{z\not\in\Lambda^\ast}\hat\pi^\ast_z({\epsilon}(t))\\
\nonumber&=\bm{\sum}\limits_{z\in\Lambda^\ast}|\hat\pi^\ast_z({\epsilon}(t))-\pi^\ast_z|+1-\bm{\sum}\limits_{z\in\Lambda^\ast}\hat\pi^\ast_z({\epsilon}(t))\\
\nonumber&=\bm{\sum}\limits_{z\in\Lambda^\ast}|\hat\pi^\ast_z({\epsilon}(t))-\pi^\ast_z|+\bm{\sum}\limits_{z\in\Lambda^\ast}\pi^\ast_z-\bm{\sum}\limits_{z\in\Lambda^\ast}\hat\pi^\ast_z({\epsilon}(t))\\
\nonumber&\le \bm{\sum}\limits_{z\in\Lambda^\ast}|\hat\pi^\ast_z({\epsilon}(t))-\pi^\ast_z|+\bm{\sum}\limits_{z\in\Lambda^\ast}|\hat\pi^\ast_z({\epsilon}(t))-\pi^\ast_z|\\
\nonumber&=2 \bm{\sum}\limits_{z\in\Lambda^\ast}\left| \frac{b^z_k{\epsilon}^c(t)^k+\cdots+b^z_h{\epsilon}^c(t)^h}{ b_k{\epsilon}^c(t)^k+\cdots+b_h{\epsilon}^c(t)^h•}-\frac{b^z_k}{b_k•} \right|\\
\nonumber&=2 \bm{\sum}\limits_{z\in\Lambda^\ast}\left| \frac{L_z({\epsilon}^c(t)^{k+1},\cdots,{\epsilon}^c(t)^{h})}{L({\epsilon}^c(t)^k,\cdots,{\epsilon}^c(t)^{h})•} \right.\\
\nonumber& \cut{\quad} \left. + \frac{{\epsilon}^c(t)^k(b^z_kb_k-b_kb^z_k)}{L({\epsilon}^c(t)^k,\cdots,{\epsilon}^c(t)^{h})•} \right|\\
&=2 {\epsilon}^c(t) \bm{\sum}\limits_{z\in\Lambda^\ast}\left| \frac{L_z(1,{\epsilon}^c(t),\cdots,{\epsilon}^c(t)^{h-k-1})}{L(1,{\epsilon}^c(t),\cdots,{\epsilon}^c(t)^{h-k})•} \right|,
\end{align}
where $L_z$ and $L$ are linear functions, the constant term of $L$ is non-zero. And by Assumption~\ref{asu6} - (1), for any $z\in \Lambda^\ast$, $L_z(1,{\epsilon}^c(t),\epsilon_a(t)^{h-k})$ and $L(1,{\epsilon}^c(t),{\epsilon}^c(t)^{h-k-1})$ converge as $t\to\infty$ because $\lim\limits_{t\to\infty}{\epsilon}^c(t)=0$. Also because $\Lambda^\ast$ contains finite elements, then \\$2\bm{\sum}\limits_{z\in\Lambda^\ast}\left| \frac{L_z(1,{\epsilon}^c(t),\cdots,{\epsilon}^c(t)^{h-k-1})}{L(1,{\epsilon}^c(t),\cdots,{\epsilon}^c(t)^{h-k})•} \right|$ is uniformly bounded. Remember that we fix $i_0$ at the beginning. And we can always chooses a constant \\$C_\epsilon\ge2\bm{\sum}\limits_{z\in\Lambda^\ast}\left| \frac{L_z(1,{\epsilon}^c(t),\cdots,{\epsilon}^c(t)^{h-k-1})}{L(1,{\epsilon}^c(t),\cdots,{\epsilon}^c(t)^{h-k})•} \right|$ such that $| |\hat\pi^\ast(t) - \pi^\ast ||\le C_\epsilon ||{\epsilon}(t)||_\infty$. It completes the proof of Claim 2. \qed\end{pf}

The following lemma shows that the sequence $\{||\hat\pi^\ast(t) -\hat\pi^\ast(t+1) ||\}_{t\ge 0}$ is summable and gives the explicit partial sums of the sequence when $t$ is large.
\begin{lemma}\label{lemma:star3}
If Assumptions~\ref{asu_error_o} - \ref{asu6} hold and ${e_i}(t)=0$ for all $t$ and all $i$, then $\bm{\sum}\limits_{\tau=0}^{+\infty}\left|\left|\hat\pi^\ast(\tau) -\hat\pi^\ast(\tau+1) \right|\right|<+\infty$. Moreover, there exists a $t^{i_0}$ such that the partial sum satisfies $\bm{\sum}\limits_{\tau=t^i}^{t}||\hat\pi^\ast(\tau) -\hat\pi^\ast(\tau+1) ||\le2| |\hat\pi^\ast(t^{i_0}+1) - \pi^\ast ||+2| |\hat\pi^\ast(t) - \pi^\ast ||, \forall t^i>t^{i_0}$ and $\forall t>t^i$.
\end{lemma}
\begin{pf}
In this paper, for any vector, we choose the $L^1$-norm, then
\begin{align}\label{L3eq1}
\nonumber &\bm{\sum}\limits_{\tau=0}^{+\infty}\left|\left|\hat\pi^\ast(\tau) -\hat\pi^\ast(\tau+1) \right|\right|\\
&=\bm{\sum}\limits_{\tau=0}^{+\infty} \bm{\sum}\limits_{z\in \mathcal{Z}}|\hat\pi^\ast_z({\epsilon}(\tau))-\hat\pi^\ast_z({\epsilon}(\tau+1))|.
\end{align}
Recall equation~(\ref{L2c1eq1}), for\cut{ $t\ge t_e$} any $z\in \mathcal{Z}$, $\hat\pi^\ast_z({\epsilon}(t))=\frac{\alpha_z({\epsilon}^c(t))}{\beta({\epsilon}^c(t))•}$ is a ratio of two polynomials of ${\epsilon}^c(t)$. Then the derivative of $\hat\pi^\ast_z({\epsilon}(t))$ is:
\begin{align*}
&\frac{\partial{\hat\pi^\ast_z ({\epsilon}(t))}}{\partial {\epsilon}^c(t) •}=\frac{1}{\beta({\epsilon}^c(t))•}\left( \frac{\partial{\alpha_z ({\epsilon}^c(t))}}{\partial {\epsilon}^c(t) •}\beta({\epsilon}^c(t)) \right. \\
&\cut{\quad} \left. - \alpha_z({\epsilon}^c(t))\frac{\partial{\beta ({\epsilon}^c(t))}}{\partial {\epsilon}^c(t) •} \right),
\end{align*}
where $\frac{\partial{\alpha_z ({\epsilon}^c(t))}}{\partial {\epsilon}^c(t) •}\beta({\epsilon}^c(t)) - \alpha_z({\epsilon}^c(t))\frac{\partial{\beta ({\epsilon}^c(t))}}{\partial {\epsilon}^c(t) •}$ is a polynomial of ${\epsilon}^c(t)$. So $\frac{\partial{\hat\pi^\ast_z ({\epsilon}^c(t))}}{\partial {\epsilon}^c(t) •}$ can be rewritten as
\begin{align*}
&\frac{\partial{\hat\pi^\ast_z ({\epsilon}(t))}}{\partial {\epsilon}^c(t) •}=\frac{c^z_l{\epsilon}^c(t)^l+c^z_{l+1}{\epsilon}^c(t)^{l+1}+\cdots+c^z_h{\epsilon}^c(t)^h}{\beta({\epsilon}^c(t))•},
\end{align*}
where $l\ge 0$, $c^z_l\neq 0$ and $\beta({\epsilon}^c(t))>0$. When ${\epsilon}^c(t)$ is sufficiently small, $c^z_l{\epsilon}^c(t)^l$ dominates the derivative. Therefore, $\exists \bar{\epsilon}^c >0$, such that the sign of $\frac{\partial{\hat\pi^\ast_z ({\epsilon}(t))}}{\partial {\epsilon}^c(t) •}$ is the sign of $c^z_l$, $\forall 0<{\epsilon}^c(t)\le \bar{\epsilon}^c $. By Assumption~\ref{asu6} - (1),  ${\epsilon}^c(t)$ strictly decreases to 0, there exists a ${t}_{i_0}$ such that ${\epsilon}^c({t}_{i_0})\le\bar{\epsilon}^c$.
We can define a partition of $\mathcal{Z}$ as follows:
\begin{align*}
&\mathcal{Z}_1=\{z\in \mathcal{Z}|\hat\pi^\ast_z ({\epsilon}(t))>\hat\pi^\ast_z ({\epsilon}(t+1)),\forall t\ge t^{i_0}\}\\
&\mathcal{Z}_2=\{z\in \mathcal{Z}|\hat\pi^\ast_z ({\epsilon}(t))<\hat\pi^\ast_z ({\epsilon}(t+1)),\forall t\ge t^{i_0}\}.
\end{align*}
Then equality~(\ref{L3eq1}) becomes:
\begin{align}\label{L3eq2}
\nonumber &\bm{\sum}\limits_{\tau=0}^{+\infty}||\hat\pi^\ast(\tau) -\hat\pi^\ast(\tau+1) ||\\
\nonumber&= \bm{\sum}\limits_{\tau=0}^{t^{i_0}}\bm{\sum} \limits_{z\in \mathcal{Z}}|\hat\pi^\ast_z({\epsilon}(\tau))-\hat\pi^\ast_z({\epsilon}(\tau+1))|\\
\nonumber& + \bm{\sum}\limits_{\tau=t^{i_0}+1}^{+\infty}\bm{\sum} \limits_{z\in \mathcal{Z}_1} (\hat\pi^\ast_z({\epsilon}(\tau))-\hat\pi^\ast_z({\epsilon}(\tau+1)) )\\
\nonumber&+ \bm{\sum}\limits_{\tau=t^{i_0}+1}^{+\infty}(1-\bm{\sum} \limits_{z\in \mathcal{Z}_1}\hat\pi^\ast_z ({\epsilon}(\tau+1) )- (1-\bm{\sum} \limits_{z\in \mathcal{Z}_1}\hat\pi^\ast_z({\epsilon}(\tau)) ) )\\
\nonumber&=\bm{\sum}\limits_{\tau=0}^{t^{i_0}}\bm{\sum} \limits_{z\in \mathcal{Z}}|\hat\pi^\ast_z({\epsilon}(\tau))-\hat\pi^\ast_z({\epsilon}(\tau+1))|\\
&+ 2\bm{\sum} \limits_{z\in \mathcal{Z}_1}(\hat\pi^\ast_z({\epsilon}(t^{i_0}+1))-\hat\pi^\ast_z(0))<+\infty.
\end{align} 
Now we consider the partial sum $\bm{\sum}\limits_{\tau=t^i}^{t}\left|\left|\hat\pi^\ast(\tau) -\hat\pi^\ast(\tau+1) \right|\right|$ when $t^i>t^{i_0}$ and $t>t^i$. By equality~(\ref{L3eq2}), we have:
\begin{align*}
\nonumber &\bm{\sum}\limits_{\tau=t^i}^{t}||\hat\pi^\ast(\tau) -\hat\pi^\ast(\tau+1) ||\\
\nonumber& \le 2\bm{\sum} \limits_{z\in \mathcal{Z}}|\hat\pi^\ast_z({\epsilon}(t^{i}+1))-\pi^\ast_z+\pi^\ast_z-\hat\pi^\ast_z({\epsilon}(t))|\\
&\le2\left| \left|\hat\pi^\ast(t^{i}+1) - \pi^\ast \right|\right|+2\left| \left|\hat\pi^\ast(t) - \pi^\ast \right|\right|.
\end{align*}\qed\end{pf}

\subsection{Proof of Theorem~\ref{The4}}\label{sec:pfthe1}
Lemma 2 shows that $\hat\pi^\ast(t) \to \pi^\ast$ whose support is $\Lambda^\ast$. Now we proceed to finish the proofs of Theorem~\ref{The4} by showing $\pi(t)\to \pi^\ast$ and quantifying its convergence rate.
\cut{
\begin{lemma}\label{lemma:the4}
Let Assumptions~\ref{asu6} - \ref{asu_error_o} hold. Then it holds that $\lim \limits_{t\to \infty}||\pi(t)-\pi^\ast||=0$ for any initial distribution $\pi(0)$. Moreover, there exists some $t^{\ast}_3$, for any $j$ such that $j(\mathcal{T}+1)\ge t^{\ast}_3$, the convergence rate could be quantified as follows:
\begin{align*}
&D(j)=O\left( \exp \left(-\bm{\sum}\limits_{\tau=t^{\ast}_3}^{(j-1)(\mathcal{T}+1)}\epsilon_d(\tau)\epsilon_a(\tau) (nm)\right) \right.\\
&\left.+4C \epsilon_\vee(t^{\ast}_3)+C \epsilon_\vee(j(\mathcal{T}+1))\right.\\
&\left.+12nm\bm{\sum}\limits_{k=t^{\ast}_3}^{\mathcal{T}+j(\mathcal{T}+1)} |{e}(k)|\right).
\end{align*}
\end{lemma}}

\begin{pf}%
Based on triangle inequality, we can get for any $t\ge 2$,
\begin{align}\label{ine_61}
&\left| \left| \pi(t) - \pi^\ast \right|\right|
\le ||\pi(t) - \hat\pi(t)||+||\hat\pi(t) - \pi^\ast||.
\end{align}
We want to prove that the two terms in the right-hand side of~(\ref{ine_61}) converge to 0 with certain rates.

{\bf Claim 3. $\lim\limits_{t\to\infty}\left| \left|\hat\pi(t) - \pi^\ast\right|\right|=0$ and there exists some $t^\vee$ such that for any $t^{\ast}_3>t^\vee$ and $t>t^{\ast}_3+1$, $\left| \left|  \hat\pi(t)   - \pi^\ast \right|\right|\le C_q\exp \left(-\bm{\sum}\limits_{\tau=t^{\ast}_3}^{t-1}\prod\limits_{i=1}^N\epsilon_i(t) |\mathcal{A}_i|\right)+4C_\epsilon||{\epsilon}(t^{\ast}_3)||_\infty+C_\epsilon||{\epsilon}(t)||_\infty$ for some constants $C_\epsilon, C_q>0$.}
\begin{pf}
Based on triangle inequality, we can get for all $t\ge 2$,
\begin{align}\label{L4eq1}
\left| \left| \hat\pi(t) - \pi^\ast \right|\right|
\le \left| \left| \hat\pi(t) - \hat\pi^\ast(t) \right|\right| + \left| \left|\hat\pi^\ast(t) - \pi^\ast \right|\right|,
\end{align}where $\hat\pi(t)$ is the distribution on $\mathcal{Z}$ at $t$ when the exploration deviations ${e_i}(t)=0$ for all $i\in\mathcal{V}$. Let $x(t)\triangleq\left|\left|\hat\pi(t) -\hat\pi^\ast(t) \right|\right| $ and $y(t)\triangleq\left|\left|\hat\pi^\ast(t) - \pi^\ast \right|\right|$. 

Let us first consider $x(t)$. Note that ${ \hat\pi^\ast(t) }^T P^{{\epsilon}(t)}={ \hat\pi^\ast(t) }^T$. Then we have:
\begin{align}\label{L4eq2}
\nonumber& x(t)=\left|\left|\hat\pi(t)-\hat\pi^\ast(t) \right|\right|\\
\nonumber&=\left|\left|\hat\pi(t)- \hat\pi^\ast(t-1) +\hat\pi^\ast(t-1) -\hat\pi^\ast(t) \right|\right|\\
\nonumber&\le \left|\left| \{ P^{{\epsilon}(t-1)}\}^T\hat\pi(t-1)-\{ P^{{\epsilon}(t-1)}\}^T\hat\pi^\ast(t-1) \right|\right|\\
 &\cut{\quad}+ \left|\left|\hat\pi^\ast(t-1) -\hat\pi^\ast(t) \right|\right|.
\end{align}
By (\ref{Lem1e1}) in the proof of Lemma~\ref{lemma:nonzero2}, \cut{when $t_q\ge t_e$, }the nonzero entries in $\{ P^{{\epsilon}(t-1)}\}^T$ can be represented as polynomials of $\{\epsilon_i(t-1), 1-\epsilon_i(t-1)\}$. Taking the nonzero entry $\prod\limits_{i=1}^N\frac{\epsilon_i(t-1) }{|\mathcal{A}_i|}$, we can decompose $\{ P^{{\epsilon}(t-1)}\}^T$ into the following:
\begin{align*}
&\{ P^{{\epsilon}(t-1)}\}^T=\prod\limits_{i=1}^N\frac{\epsilon_i(t-1) }{|\mathcal{A}_i|}Q+R(t-1),
\end{align*}where $Q$ is a $|\mathcal{Z}|\times |\mathcal{Z}|$ matrix with all entries are 1. Because $ P^{{\epsilon}(t-1)}$ is a transition matrix, then the $\{ P^{{\epsilon}(t-1)}\}^T$ is a column stochastic matrix where each column sum is equal to 1. It follows that the column sums of $\prod\limits_{i=1}^N\frac{\epsilon_i(t-1) }{|\mathcal{A}_i|}Q$ equal $\prod\limits_{i=1}^N\frac{\epsilon_i(t-1) }{|\mathcal{A}_i|}|\mathcal{Z}|=\prod\limits_{i=1}^N\epsilon_i(t-1)|\mathcal{A}_i|$ since $|\mathcal{Z}|=\left(\prod\limits_{i=1}^N|\mathcal{A}_i|\right)^2$, and the column sums of $R(t-1)$ equal $c(t-1)=1-\prod\limits_{i=1}^N\epsilon_i(t-1)|\mathcal{A}_i|$.

By (1) in Assumption~\ref{asu6}, $\prod\limits_{i=1}^N\epsilon_i(t-1)|\mathcal{A}_i|$ strictly decreases to 0. Then there exists a $t^{|\mathcal{A}|}$ such that $\prod\limits_{i=1}^N\epsilon_i(t-1)|\mathcal{A}_i|<1$ for all $t\ge  t^{|\mathcal{A}|}$, which implies $0<c(t-1)<1$ for all $t\ge t^{|\mathcal{A}|}$ and the column sums of $c(t-1)^{-1}R(t-1)$ equal 1.
%
Let $v(t-1)\triangleq \left|\left|\hat\pi^\ast(t-1) -\hat\pi^\ast(t) \right|\right|$. And consider $t\ge t^{|\mathcal{A}|}$, then inequality~(\ref{L4eq2}) becomes:
\begin{align}\label{L4eq8}
\nonumber &x(t)\le  ||\{P^{{\epsilon}(t-1)}\}^T(\hat\pi(t-1)- \hat\pi^\ast(t-1)) ||+ v(t-1)\\
\nonumber &=||\prod\limits_{i=1}^N\frac{\epsilon_i(t-1) }{|\mathcal{A}_i|}Q(\hat\pi(t-1)- \hat\pi^\ast(t-1))\\
 & \cut{\quad}+R(t-1)(\hat\pi(t-1)-\hat\pi^\ast(t-1))||+ v(t-1),
\end{align}where $\hat\pi(t-1)$ and $\hat\pi^\ast(t-1)$ are both stochastic vectors whose sum of elements is equal to one. And by the construction of $Q$, we have $Q(\hat\pi(t-1)- \hat\pi^\ast(t-1))=0$. Then inequality~(\ref{L4eq8}) becomes:
\begin{align*}
&x(t)\le c(t-1)||c(t-1)^{-1}R(t-1)\\
&\cut{\quad}\times(\hat\pi(t-1)- \hat\pi^\ast(t-1))||+ v(t-1)\\
 &\le c(t-1)||c(t-1)^{-1}R(t-1)||\\
 &\cut{\quad}\times||\hat\pi(t-1)- \hat\pi^\ast(t-1)||+ v(t-1)\\
 &= c(t-1)x(t-1)+ v(t-1),
\end{align*}where $c(t-1)\in (0,1)$ for all $ t\ge t^{|\mathcal{A}|}$. With inequality $\log(1-x)<-x, \forall x\in (0,1)$, for any $t>t^{|\mathcal{A}|}$ and $t^{\ast}_3\ge t^{|\mathcal{A}|}$, we have :
\begin{align}\label{L4eq3}
\nonumber &x(t)\le  c(t-1)x(t-1) +v(t-1)\\
\nonumber &\le \prod\limits_{\tau=t^{\ast}_3}^{t-1}c(\tau) x(t^{\ast}_3)+v(t-1)+\bm{\sum}\limits_{\tau=t^{\ast}_3}^{t-2}\left( \prod\limits_{i=\tau+1}^{t-1}c(i)v(\tau) \right)\\
\nonumber &\le x(t^{\ast}_3)\prod\limits_{\tau=t^{\ast}_3}^{t-1}\exp\left(-(1-c(\tau))\right)+v(t-1)\\
\nonumber &\cut{\quad}+\bm{\sum}\limits_{\tau=t^{\ast}_3}^{t-2}\left( \prod\limits_{i=\tau+1}^{t-1} \exp\left(-(1-c(i))\right) v(\tau) \right)\\
&\le x(t^{\ast}_3)\exp\left(-\bm{\sum}\limits_{\tau=t^{\ast}_3}^{t-1}(1-c(\tau))\right)+\bm{\sum}\limits_{\tau=t^{\ast}_3}^{t-1}v(\tau).
\end{align}%
Note that inequality (\ref{L4eq3}) holds for any $t^{\ast}_3\ge t^{|\mathcal{A}|}$. And by Lemma~\ref{lemma:star3}, $v(\tau)$ is summable, we first take the limit of $t$ and then take the limit of $t^{\ast}_3$, by the summability of $v(\tau)$, we can have $\lim\limits_{t^{\ast}_3\to \infty}\lim\limits_{t\to \infty}\bm{\sum}\limits_{\tau=t^{\ast}_3}^{t-1}v(\tau)=0$. And for any $t^i>t^{i_0}$ and $t> t^i$, $\bm{\sum}\limits_{\tau=t^i}^{t}||\hat\pi^\ast(\tau) -\hat\pi^\ast(\tau+1) ||\le  2| |\hat\pi^\ast(t^{i}+1) - \pi^\ast ||+2| |\hat\pi^\ast(t) - \pi^\ast ||$.

Let $t^\vee\triangleq\max\{t^{i_0},t^{|\mathcal{A}|}\}$. Since Inequality (\ref{L4eq3}) holds for any $t^{\ast}_3\ge t^{|\mathcal{A}|}$, we can take $t^{\ast}_3>t^\vee$. Recall that $1-c(\tau)=\prod\limits_{i=1}^N\epsilon_i(\tau)|\mathcal{A}_i|$, then for any $t>t^{\ast}_3$, inequality (\ref{L4eq3}) becomes:
%
\begin{align}\label{L4eq6}
\nonumber &x(t)\le x(t^{\ast}_3)  \exp \left(-\bm{\sum}\limits_{\tau=t^{\ast}_3}^{t-1}\prod\limits_{i=1}^N\epsilon_i(\tau)|\mathcal{A}_i|\right)\\
\nonumber &\cut{\quad} + 2| |\hat\pi^\ast(t^{\ast}_3+1) - \pi^\ast ||+2| |\hat\pi^\ast(t-1) - \pi^\ast ||\\
\nonumber&= x(t^{\ast}_3) \exp \left(-\bm{\sum}\limits_{\tau=t^{\ast}_3}^{t-1}\prod\limits_{i=1}^N\epsilon_i(\tau)|\mathcal{A}_i|\right)\\
&\cut{\quad} +2y(t^{\ast}_3+1)+2y(t-1).
\end{align}%
Combining inequalities (\ref{L4eq1}) and (\ref{L4eq6}), we can get for any $t^{\ast}_3>t^\vee$:
\begin{align}\label{L4eq7}
\nonumber&\left| \left| \hat\pi(t) - \pi^\ast \right|\right|\le x(t^{\ast}_3)\exp \left(-\bm{\sum}\limits_{\tau=t^{\ast}_3}^{t-1}\prod\limits_{i=1}^N\epsilon_i(\tau)|\mathcal{A}_i|\right)\\
&\cut{\quad} +2y(t^{\ast}_3+1)+2y(t-1)+y(t).
\end{align}
By (2) in Assumption~\ref{asu6}, $\prod\limits_{i=1}^N\epsilon_i(\tau)|\mathcal{A}_i|$ is not summable. Therefore, for any $t^{\ast}_3>t^\vee$, we have \begin{align*}
\lim\limits_{t\to\infty}x(t^{\ast}_3)\exp\left(-\bm{\sum}\limits_{\tau=t^{\ast}_3}^{t-1}\prod\limits_{i=1}^N\epsilon_i(\tau)|\mathcal{A}_i|\right)=0.
\end{align*}
By Lemma~\ref{lemma:star2}, we have
\begin{align*}
&2y(t^{\ast}_3+1)+2y(t-1)+y(t)\\
&\le 2C_\epsilon ||{\epsilon}(t^{\ast}_3+1)||_\infty+2C_\epsilon |{\epsilon}(t-1)_\infty+C_\epsilon ||{\epsilon}(t)||_\infty\\
&\le 4C_\epsilon ||{\epsilon}(t^{\ast}_3)||_\infty+C_\epsilon ||{\epsilon}(t)||_\infty,
\end{align*} where $2C_\epsilon ||{\epsilon}(t^{\ast}_3)+1||_\infty+2C_\epsilon ||{\epsilon}(t-1)||_\infty\le 4C_\epsilon ||{\epsilon}(t^{\ast}_3)||_\infty$ since $\epsilon_i(t)$ is strictly decreasing to 0. And there exists a positive constant $C_q$ such that $x(t^{\ast}_3) \le C_q$. Therefore, for any $t^{\ast}_3>t^\vee$ and $t>t^{\ast}_3+1$, (\ref{L4eq7}) becomes:
\begin{align*}
&\left| \left| \hat\pi(t) - \pi^\ast \right|\right|\le C_q\exp \left(-\bm{\sum}\limits_{\tau=t^{\ast}_3}^{t-1}\prod\limits_{i=1}^N\epsilon_i(\tau)|\mathcal{A}_i|\right)\\
&\cut{\quad} +4C_\epsilon ||{\epsilon}(t^{\ast}_3)||_\infty+C_\epsilon ||{\epsilon}(t)||_\infty.
\end{align*}
Therefore we reach Claim 3.
\cut{
Therefore we reach that $\left| \left| \hat\pi(t_q) - \pi^\ast \right|\right|$ converges to 0 as $t\to\infty$ and for any $t>t_q^\vee+\mathcal{T}+1$ and $t^{\ast}_3>t_q^\vee$, \begin{align*}
&\left| \left|\hat\pi(t_q) - \pi^\ast \right|\right|\\
&= O\left(\epsilon_\vee(t_q^\ast)+\epsilon_\vee(t_q)+\exp \left(-\bm{\sum}\limits_{\tau_q=t^{\ast}_3}^{t-\mathcal{T}-1}\epsilon_d(\tau_q)\epsilon_a(\tau_q) (nm)\right)\right).
\end{align*}}
\qed\end{pf}

{\bf Claim 4. $\lim\limits_{t\to\infty}| | \pi(t)  - \hat\pi(t)  ||=0$ and there exists some $t^c\cut{\ge t_e+\mathcal{T}+1}$ such that for any $t^{\ast}_4\ge t^c+1$ and $t\ge t^{\ast}_4$, 
$|| \pi(t)  - \hat\pi(t)  ||\le C_c \exp (- \bm{\sum}\limits_{\tau=t^{\ast}_4}^{t-1}\prod\limits_{i=1}^N\tilde\epsilon_i(\tau)|\mathcal{A}_i|)+4^N e_r(t^{\ast}_4)$ for some constant $C_c>0$.}
\begin{pf}
Based on triangle inequality, we can get for all $t\ge 2$,
\begin{align}\label{c4e1}
\nonumber&\left| \left|\pi(t)- \hat\pi(t)  \right|\right|
\le \left| \left| \pi(t)- \{ P^{\tilde{\epsilon}(t-1)}\}^T\hat\pi(t-1) \right|\right| \\
&\cut{\quad}+ \left| \left| \{ P^{\tilde{\epsilon}(t-1)}\}^T\hat\pi(t-1) - \hat\pi(t)\right|\right|.
\end{align}
With $\pi(t)=\{ P^{\tilde{\epsilon}(t-1)}\}^T\pi(t-1)$ and $\hat\pi(t)=\{P^{{\epsilon}(t-1)}\}^T\allowbreak\hat\pi(t-1)$, (\ref{c4e1}) becomes:
\begin{align}\label{c4e2}
\nonumber&\left| \left|\pi(t)- \hat\pi(t)  \right|\right|\le \left|\left| \{ P^{\tilde{\epsilon}(t-1)}\}^T\left(\pi(t-1)-\hat\pi(t-1)\right) \right|\right|\\
 &+ \left|\left|\left(\{ P^{\tilde{\epsilon}(t-1)}\}^T-\{ P^{{\epsilon}(t-1)}\}^T\right)\hat\pi(t-1) \right|\right|.
\end{align}
Based on Lemma~\ref{lemma:nonzero2},\cut{ when $t_q\ge t_e$} the nonzero entries in $\{ P^{\tilde{\epsilon}(t-1)}\}^T$ can be represented as polynomials of $\{\tilde\epsilon_i(t-1), 1-\tilde\epsilon_i(t-1)\}$. Taking the nonzero entry $\prod\limits_{i=1}^N\frac{\tilde\epsilon_i(t-1) }{|\mathcal{A}_i|}$, we can decompose $\{ P^{\tilde{\epsilon}(t-1)}\}^T$ into the following:
$
\{ P^{\tilde{\epsilon}(t-1)}\}^T=\prod\limits_{i=1}^N\frac{\tilde\epsilon_i(t-1)}{|\mathcal{A}_i|}Q+R'(t-1)$, where $Q$ is a $|\mathcal{Z}|\times |\mathcal{Z}|$ matrix with all entries are 1. Because $ P^{\tilde{\epsilon}(t-1)}$ is a transition matrix, then the $\{ P^{\tilde{\epsilon}(t-1)}\}^T$ is a column stochastic matrix where each column sum is equal to 1. It follows that the column sums of $\prod\limits_{i=1}^N\frac{\tilde\epsilon_i(t-1) }{|\mathcal{A}_i|}Q$ equal $\prod\limits_{i=1}^N\frac{\tilde\epsilon_i(t-1) }{|\mathcal{A}_i|}|\mathcal{Z}|=\prod\limits_{i=1}^N\tilde\epsilon_i(t-1) |\mathcal{A}_i|$, and the column sums of $R'(t-1)$ equal $c(t-1)=1-\prod\limits_{i=1}^N\tilde\epsilon_i(t-1)|\mathcal{A}_i|$.

Let $x(t)=\left| \left|\pi(t)- \hat\pi(t)  \right|\right|$. And by the structure of $Q$ and the fact that $\pi(t-1)$ and $\hat\pi(t-1)$ are both stochastic vectors whose sum of elements is equal to 1, we have $Q(\pi(t-1)- \hat\pi(t-1))=0$. Then (\ref{c4e2}) becomes:
\begin{align}\label{c4e3}
\nonumber&x(t)\le||R'(t-1)(\pi(t-1)-\hat\pi(t-1))||\\
 \nonumber&\cut{\quad}+\left|\left|\left(\{ P^{\tilde{\epsilon}(t-1)}\}^T-\{ P^{{\epsilon}(t-1)}\}^T\right)\hat\pi(t-1) \right|\right|\\
\cut{\nonumber &\le c(t-1)x(t-1)\\
\nonumber&  +\left|\left|\left(\{ P^{\tilde{\epsilon}(t-1)}\}^T-\{ P^{{\epsilon}(t-1)}\}^T\right)\hat\pi(t-1) \right|\right|.\\}
 &\le c(t-1)x(t-1)+\left|\left|\{ P^{\tilde{\epsilon}(t-1)}\}^T-\{ P^{{\epsilon}(t-1)}\}^T \right|\right|.
\end{align}
By (\ref{Lem1e1}) in the proof of Lemma~\ref{lemma:nonzero2},\cut{ for any $t-1\ge t_e$,} any entry in $ P^{\tilde{\epsilon}(t-1)}$ can be represented as a summation of at most $\cut{\bm{\sum}\limits_{j=0}^N {N\choose j}=}2^N$ polynomials\cut{ (this happens when $\mathcal{V}_{er}(x,y)=\emptyset$)}. And each polynomial is a product of $N$ monomials; e.g., $\prod\limits_{i=1}^N\frac{\tilde\epsilon_i(t-1)}{|\mathcal{A}_i|}$. And the entries in  $P^{{\epsilon}(t-1)}$ have the same form with $\tilde\epsilon_i(t-1)=\epsilon_i(t-1)$. Then the difference of any pair of entries $( P^{\tilde{\epsilon}(t-1)}(x,y), P^{{\epsilon}(t-1)}(x,y))$ has at most $4^{N}$ terms (for example, $\prod\limits_{i=1}^N\frac{\tilde\epsilon_i(t-1)}{|\mathcal{A}_i|}-\prod\limits_{i=1}^N\frac{\epsilon_i(t-1)}{|\mathcal{A}_i|}$ has $2^N-1$ terms). And each term is less than $\prod\limits_{i=1}^N\frac{||{e}(t-1)||_\infty}{|\mathcal{A}_i|}$, where $||{e}(t)||_\infty = \max\{|e^1_c(t)|,\cdots, |e^N_c(t)|\}$. Then any pair of entries $( P^{\tilde{\epsilon}(t-1)}(x,y), P^{{\epsilon}(t-1)}(x,y))$ satisfy that:
\begin{align}
\nonumber|  P^{\tilde{\epsilon}(t-1)}(x,y)-P^{{\epsilon}(t-1)}(x,y) |\le 4^N \prod\limits_{i=1}^N\frac{1}{|\mathcal{A}_i|} ||{e}(t-1)||_\infty.
\end{align}
\cut{Therefore, $||\{ P^{\tilde{\epsilon}(t-1)}\}^T-\{ P^{{\epsilon}(t-1)}\}^T ||\le \max\limits_{x}\bm{\sum}\limits_{y} |  P^{\tilde{\epsilon}(t-1)}(x,y)-P^{{\epsilon}(t-1)}(x,y) |\le  4^N||{e}(t-1)||_\infty^N\prod\limits_{i=1}^N|\mathcal{A}_i| $.} Then (\ref{c4e3}) becomes:
\begin{align}\label{c4e5}
x(t)\le c(t-1)x(t-1) + 4^N||{e}(t-1)||_\infty^N\prod\limits_{i=1}^N|\mathcal{A}_i|.
\end{align}
Let $\chi(t)=x(t)-\frac{ 4^N||{e}(t)||_\infty^N\prod\limits_{i=1}^N|\mathcal{A}_i|}{1- c(t)}=x(t)-4^N e_r(t)\cut{ =x(t)-\frac{2^N (2^N-1)||{e}(t)||_\infty^N}{\prod\limits_{i=1}^N(\epsilon_i(t)+{e_i}(t) )} }$. Then from (\ref{c4e5}), we can get:
\begin{align}\label{c4e6}
\chi(t)\le 
c(t-1) \chi(t-1)+4^N e_r(t-1)-4^N e_r(t). 
\end{align}

Inequality (\ref{c4e6}) holds for any $t\ge 2$. Recall that $\tilde\epsilon_i(t)\in(0,1],\forall i\in \mathcal{V},\forall t$, then $c(t-1)\le 1$. By simple algebraic operations, we have $c(t-1) \ge 1-\prod\limits_{i=1}^N(\epsilon_i(t-1)+||{e}(t-1)||_\infty)|\mathcal{A}_i|$. And by (1) and (3) in Assumption~\ref{asu6}, $\prod\limits_{i=1}^N(\epsilon_i(t-1)+||{e}(t-1)||_\infty)$ converges to 0. Therefore, there exists a $t^c$ such that $\prod\limits_{i=1}^N(\epsilon_i(t-1)+||{e}(t-1)||_\infty)|\mathcal{A}_i|<1,\forall t\ge t^c+1$. Then $c(t-1)\in (0,1), \forall t\ge t^c+1$. By manipulating inequality $\log(1-x)<-x, \forall x\in (0,1)$ and $\exp(-x)<1, \forall x>0$, (\ref{c4e6})  can be rewritten for all $t^{\ast}_4\ge t^c+1$ and $t\ge t^{\ast}_4$ as follows:
\begin{align}
\nonumber&\chi(t)\le \prod\limits_{\tau=t^{\ast}_4}^{t-1} c(\tau) \chi(t^{\ast}_4)+4^N e_r(t-1)-4^N e_r(t)\\
\nonumber& \cut{\quad}+ \bm{\sum}\limits_{\tau=t^{\ast}_4}^{t-2}\left( \prod\limits_{j=\tau+1}^{t-1}c(j) \right)4^N\left( e_r(\tau)-e_r(\tau+1)\right)\\
\nonumber&\le  \chi(t^{\ast}_4)\exp \left(- \bm{\sum}\limits_{\tau=t^{\ast}_4}^{t-1}(1- c(\tau))\right)\\
\nonumber&\cut{\quad} + \bm{\sum}\limits_{\tau=t^{\ast}_4}^{t-1}4^N\left( e_r(\tau)- e_r(\tau+1)\right)\\
\nonumber&\le \chi(t^{\ast}_4)\exp \left(- \bm{\sum}\limits_{\tau=t^{\ast}_4}^{t-1}(1- c(\tau))\right) +4^N(e_r(t^{\ast}_4)-e_r(t)). 
\end{align}%
Plug $\chi(t)=x(t)-4^N e_r(t)$ and $c(t)=1-\prod\limits_{i=1}^N\tilde\epsilon_i(t)|\mathcal{A}_i|$ in the above inequality, and we have:
\begin{align}\label{c4e8}
&x(t)\le \chi(t^{\ast}_4)\exp \left(- \bm{\sum}\limits_{\tau=t^{\ast}_4}^{t-1}\prod\limits_{i=1}^N\tilde\epsilon_i(\tau)|\mathcal{A}_i|\right) +4^N e_r(t^{\ast}_4). 
\end{align}%
\cut{The second inequality holds because ${e_i}(\tau)$ and $\frac{ 2^N (2^N-1)(||{e}(t)||_\infty)^N}{\prod\limits_{i=1}^N\epsilon_i(t)}$ are both non-negative strictly decreasing (Assumption~\ref{asu_error_c_new}). }%
By Assumption~\ref{asu6} - (2),  $\prod\limits_{i=1}^N\tilde\epsilon_i(\tau)|\mathcal{A}_i|$ is not summable. Therefore,
$\lim\limits_{t\to\infty}\chi(t^{\ast}_4)\exp\left(-\bm{\sum}\limits_{\tau=t^{\ast}_4}^{t-1}\prod\limits_{i=1}^N\tilde\epsilon_i(\tau)|\mathcal{A}_i|\right)=0$.
And by Assumption~\ref{asu6} - (3), $\lim\limits_{t^{\ast}_4\to\infty}4^N e_r(t^{\ast}_4)=0$.  We first take the limit of $t$ and then take the limit of $t^{\ast}_4$, we can have $\lim\limits_{t^{\ast}_4\to\infty}\lim\limits_{t\to\infty}\chi(t^{\ast}_4)\exp\left(-\bm{\sum}\limits_{\tau=t^{\ast}_4}^{t-1}\prod\limits_{i=1}^N\tilde\epsilon_i(\tau)|\mathcal{A}_i|\right)+4^N e_r(t^{\ast}_4)=0$.
And there exists a positive constant $C_c$ such that $\chi(t^{\ast}_4) \le C_c$. Then for $t^{\ast}_4\ge t^c+1$ and $t\ge t^{\ast}_4$:
\begin{align*}
&|| \pi(t)  - \hat\pi(t)  ||\le C_c \exp (- \bm{\sum}\limits_{\tau=t^{\ast}_4}^{t-1}\prod\limits_{i=1}^N\tilde\epsilon_i(\tau)|\mathcal{A}_i|)+4^N e_r(t^{\ast}_4). 
\end{align*}
Therefore we reach Claim 4. \qed
\end{pf}
Combining Claim 3 and Claim 4, we get that for Markov chain $\mathcal{M}$, its state distribution $\{\pi(t)\}$ converges to limiting distribution $\pi^\ast$. Moreover, by triangle inequality, Claim 3 and Claim 4, there exists some $t_{min}=\max \{t^\vee,t^c\}$ and ${C}\ge\max\{C_q,4C_\epsilon,C_c,4^N\}$, such that for any $t^\ast>t_{min}$ and $t>t^\ast+1$, inequality~(\ref{theorem_in}) holds. 
It completes the proof of Theorem~\ref{The4}. \qed\end{pf}

\subsection{Proof of Corollary~\ref{Cor1}}\label{sec:pfCor1}
\begin{pf}
From the last two paragraph of Section~\ref{sec:pfthe1}, the constant $C$ in inequality~(\ref{theorem_in}) can be estimated as ${C}\ge\max\{C_q,4C_\epsilon,C_c,4^N\}$. Now we will prove that there exists a set of feasible constants $C_q,4C_\epsilon,C_c$ such that $\max\{|\mathcal{Z}|,4^N, 32N|\mathcal{Z}|^{|\mathcal{Z}|+4}(N+1)^{|\mathcal{Z}|}2^{N|\mathcal{Z}|}\frac{C_{max} }{C_{min} }\}\ge \max\{C_q,4C_\epsilon,C_c,4^N\}$.

From Claim 3, $C_q$ can be any constant that satisfies $C_q\ge \left|\left|\hat\pi(t^{\ast}_3)-\hat\pi^\ast(t^{\ast}_3) \right|\right|=\bm{\sum}\limits_{z\in \mathcal{Z}}|\hat\pi_z({\epsilon}(t^{\ast}_3))-\hat\pi^\ast_z({\epsilon}(t^{\ast}_3))|$. Since $\hat\pi_z({\epsilon}(t^{\ast}_3))\in[0,1]$ and $\hat\pi^\ast_z({\epsilon}(t^{\ast}_3))\in[0,1]$. $C_q$ can be $C_q=|\mathcal{Z}|\ge \left|\left|\hat\pi(t^{\ast}_3)-\hat\pi^\ast(t^{\ast}_3) \right|\right|$. And from Claim 4, $C_c$ can be any constant that satisfies $C_c\ge\chi(t^{\ast}_4) = \left| \left|\pi(t^{\ast}_4)- \hat\pi(t^{\ast}_4)  \right|\right|-4^N e_r(t^{\ast}_4)$. Here, $\left| \left|\pi(t^{\ast}_4)- \hat\pi(t^{\ast}_4)  \right|\right|-4^N e_r(t^{\ast}_4)\le \left| \left|\pi(t^{\ast}_4)- \hat\pi(t^{\ast}_4)  \right|\right|\le \bm{\sum}\limits_{z\in \mathcal{Z}}|\pi_z({\epsilon}(t^{\ast}_4))-\hat\pi_z({\epsilon}(t^{\ast}_4))|\le |\mathcal{Z}|$. Then $C_c$  can be  $C_c= |\mathcal{Z}|$. 

From Claim 2, $C_\epsilon$ can be any constant that satisfies $C_\epsilon\ge 2\bm{\sum}\limits_{z\in\Lambda^\ast}| \frac{L_z(1,{\epsilon}^c(t),\cdots,{\epsilon}^c(t)^{h-k-1})}{L(1,{\epsilon}^c(t),\cdots,{\epsilon}^c(t)^{h-k})•} |=\underline{C_\epsilon}$. \cut{For presentation simplicity, denote by $\underline{C_\epsilon} \triangleq 2\bm{\sum}\limits_{z\in\Lambda^\ast}| \frac{L_z(1,{\epsilon}^c(t),\cdots,{\epsilon}^c(t)^{h-k-1})}{L(1,{\epsilon}^c(t),\cdots,{\epsilon}^c(t)^{h-k})•} |$ }And by (\ref{L2c2eq1}), we have:
\begin{align}\label{Cor1in1}
\nonumber&\underline{C_\epsilon}=2 \bm{\sum}\limits_{z\in\Lambda^\ast}\frac{1}{{\epsilon}^c(t)}\left| \frac{b^z_k{\epsilon}^c(t)^k+\cdots+b^z_h{\epsilon}^c(t)^h}{ b_k{\epsilon}^c(t)^k+\cdots+b_h{\epsilon}^c(t)^h•}-\frac{b^z_k}{b_k•} \right|\\
\nonumber&=2 \bm{\sum}\limits_{z\in\Lambda^\ast}\left| \frac{b_k(b^z_{k+1}{\epsilon}^c(t)^{k}+\cdots+b^z_h{\epsilon}^c(t)^{h-1})}{ b_k(b_k{\epsilon}^c(t)^k+\cdots+b_h{\epsilon}^c(t)^h)} \right.\\
\nonumber& \cut{\quad} \left. -\frac{b^z_k(b_{k+1}{\epsilon}^c(t)^{k}+\cdots+b_h{\epsilon}^c(t)^{h-1})}{b_k(b_k{\epsilon}^c(t)^k+\cdots+b_h{\epsilon}^c(t)^h)}\right|\\
%
\nonumber&\le 2 \bm{\sum}\limits_{z\in\mathcal{Z}} \left(\left| \frac{|b^z_{k+1}|+\cdots+|b^z_h|}{ b_k+\cdots+b_h{\epsilon}^c(t)^{h-k}}  \right|\right.\\
& \cut{\quad} \left. +\left|\frac{|b_{k+1}|+\cdots+|b_h|}{b_k+\cdots+b_h{\epsilon}^c(t)^{h-k}}\right|\right),
\end{align} 
where we use $|{\epsilon}^c(t)\le|\le 1$ in the inequality.
Based on equalities (\ref{S4e1}) and (\ref{L2c1eq1}), we have the following relation:
\begin{align}\label{Cor1e1}
\nonumber&b^z_k{\epsilon}^c(t)^k+\cdots+b^z_h{\epsilon}^c(t)^h=\sigma_z({\epsilon}(t))\\
&=\bm{\sum}\limits_{T\in G_{{\epsilon}(t)}(z)}\prod\limits_{(x,y)\in E(T)}P^{{\epsilon}(t)}(x,y).
\end{align}
Here we assume $w_i(t)=0, \forall i\in\mathcal{V} t\ge 0$. Then the mis-exploitations happen with probability 0. So 
\begin{align}\label{Cor1e2}
\nonumber&P^{\epsilon(t)}(x,y)=\prod\limits_{ i\in \mathcal{V}_{er}(x,y)}\frac{\gamma_i{\epsilon}^c(t)}{|\mathcal{A}_i|}\\
&\cut{\quad}\times\prod\limits_{j\in \mathcal{V}_{ex}(x,y)}(1+(-\gamma_j+ \frac{\gamma_j}{|\mathcal{A}_j|}){\epsilon}^c(t)).
\end{align}
So $P^{{\epsilon}(t)}(x,y)$ is a polynomial of ${\epsilon}^c(t)$ and can be written as $P^{{\epsilon}(t)}(x,y)=d_0(x,y)+d_{1}(x,y){\epsilon}^c(t)+\cdots+d_{N}(x,y){\epsilon}^c(t)^N$. Note that some coefficient $d_{m}(x,y)$ could be $0$. Denote by $\hat d(x,y)$ the coefficient of the least degree term in $P^{{\epsilon}(t)}(x,y)$. In fact, coefficients $b^z_{k},\cdots b^z_{h}, b_k,\cdots, b_h$ consists of $d_{m}(x,y)$, where $m\in\{0,\cdots,N\}$. In Claim 5, we will first find an upper bound of $|d_{m}(x,y)|$ and a lower bound of $|\hat d(x,y)|$. In Claims 6 and 7, we will estimate the last sum in inequality (\ref{Cor1in1}) by finding an upper bound of  $|b^z_{k+1}|+\cdots+|b^z_h|$ and $|b_{k+1}|+\cdots+|b_h|$ and a lower bound of $|b_k+\cdots+b_h{\epsilon}^c(t)^{h-k}|$.


{\bf Claim 5. For any $x,y\in\mathcal{Z}$ and $m\in \{0,\cdots,N\}$, $|d_m(x,y)|\le 2^N\max\{1,||\gamma||_\infty^N\}$. And for any $x,y\in\mathcal{Z}$, $\hat d(x,y)\ge \min\{(\frac{|\gamma|_{min}}{|\mathcal{A}|_{\infty}})^N, 1\}$.}

\begin{pf}
\cut{ The coefficients of $P^{{\epsilon}(t)}(x,y)$; e.g., $d_l(x,y)$, are linear combinations of $\frac{\gamma_i}{|\mathcal{A}_i|}$, $1$, $-\gamma_j + \frac{\gamma_j}{|\mathcal{A}_j|}$ and their products.} \cut{By binomial theorem~\cite{abramowitz1964handbook}, }Expanding the right hand side of (\ref{Cor1e2}) yields the sum of at most $2^N$ monomials where each monomial is a product of $\prod\limits_{ i\in \mathcal{V}_{er}(x,y)}\frac{\gamma_i{\epsilon}^c(t)}{|\mathcal{A}_i|} $, $1$ and $(-\gamma_j+ \frac{\gamma_j}{|\mathcal{A}_j|}){\epsilon}^c(t)$. 
\begin{align*}
&d_m(x,y)=\bm{\sum}\limits_{ \{\mathcal{V}'(x,y)\subseteq \mathcal{V}_{ex}(x,y)| | \mathcal{V}_{er}(x,y)|+| \mathcal{V}'(x,y)|=m\} }\\
&\left(\prod\limits_{ i\in \mathcal{V}_{er}(x,y)}\frac{\gamma_i}{|\mathcal{A}_i|} \prod\limits_{j\in \mathcal{V}'(x,y)}(-\gamma_{j}+ \frac{\gamma_{j}}{|\mathcal{A}_{j}|})\right),
\end{align*}And there are $ {|\mathcal{V}_{ex}(x,y)|}\choose {m-|\mathcal{V}_{er}(x,y)|}$ choices of $\mathcal{V}'(x,y)$. Since $|\mathcal{V}_{ex}(x,y)|\le N$ here we use the upper bound $ {{|\mathcal{V}_{ex}(x,y)|}\choose {m-|\mathcal{V}_{er}(x,y)|}}\le  {{N}\choose {m-|\mathcal{V}_{er}(x,y)|}}\le 2^N$ for presentation simplicity.
%
%
%
Note that $||\gamma||_\infty = \max\limits_{i\in \mathcal{V}}\gamma_i$, then $\frac{\gamma_i}{|\mathcal{A}_i|}< ||\gamma||_\infty$ and $|-\gamma_j + \frac{\gamma_j}{|\mathcal{A}_j|}|\le ||\gamma||_\infty$.
Then for all $m\in\{0,\cdots,N\}$, 
\begin{align}
\nonumber&|d_m(x,y)|\le\bm{\sum}\limits_{| \mathcal{V}_{er}(x,y)|+| \mathcal{V}'_{ex}(x,y)|=m}(\max\{1,||\gamma||_\infty\})^m\\
\nonumber&\le2^N\max\{1,||\gamma||_\infty^N\}.
\end{align}%
By equality (\ref{Cor1e2}), $\hat d(x,y)=\prod\limits_{ i\in \mathcal{V}_{er}(x,y)}\frac{\gamma_i}{|\mathcal{A}_i|}$ when $\mathcal{V}_{er}(x,y) \neq \emptyset$ or $\hat d(x,y)=1$ when $\mathcal{V}_{er}(x,y) = \emptyset$. Note that, $\hat d(x,y)>0$ for any $x,y\in\mathcal{Z}$. With $\frac{\gamma_i}{|\mathcal{A}_i|}\ge \frac{|\gamma|_{min}}{|\mathcal{A}|_{\infty}}$, $\hat d(x,y)\ge  \min\{(\frac{|\gamma|_{min}}{|\mathcal{A}|_{\infty}})^N, 1\}$.
\qed\end{pf}
%
%

{\bf Claim 6. $|b^z_{k+1}|+\cdots+|b^z_h|$ and $|b_{k+1}|+\cdots+|b_h|$ are both upper bounded by $ N|\mathcal{Z}|^{|\mathcal{Z}|+3}(N+1)^{|\mathcal{Z}|}2^{N|\mathcal{Z}|}C_{max}\cut{\max\{1,||\gamma||_\infty^{N|\mathcal{Z}|}\}}$.}

\begin{pf}
Notice that 
\begin{align}\label{Cor1e4}
\nonumber&\prod\limits_{(x,y)\in E(T)}P^{{\epsilon}(t)}(x,y)\\
\nonumber&=\prod\limits_{(x,y)\in E(T)}(d_0(x,y)+\cdots+d_{N}(x,y){\epsilon}^c(t)^N)\\
&=f^T_0+f^T_1{\epsilon}^c(t)+\cdots+f^T_{N|T|} {\epsilon}^c(t)^{N|T|}, 
\end{align}
where $|T|$ is the number of edges of tree $T$\cut{ and $f^T_m$ is the coefficient of the degree $i$ term, which could be 0}. 
For the analytical simplicity, denote the enumeration of edges in $T$ as $E(T)=\{g_1,g_2,\cdots,g_{|T|}\}$ and denote by $d_{l_g}(g)$ the coefficient of the $l_g$-th degree term in the polynomial $P^{{\epsilon}(t)}(x,y)$, where $l_g\in\{0,\cdots,N\}$. Then $\prod\limits_{(x,y)\in E(T)}(d_0(x,y)+\cdots+d_{N}(x,y){\epsilon}^c(t)^N)=\prod\limits_{g=g_1}^{g_{|T|}}(d_0(g)+\cdots+d_{N}(g){\epsilon}^c(t)^N)$, which can be expanded as the sum of $(N+1)^{|T|}$ monomials where each monomial is in the form of $\prod\limits_{g=g_1}^{g_{|T|}}d_{l_g}(g){\epsilon}^c(t)^{l_g}$. Then $f^T_m=\bm{\sum}\limits_{\{l_{g_1},\cdots,l_{g_{|T|}} | l_{g_1}+\cdots+l_{g_{|T|}}=m\}} \left(\prod\limits_{g=g_1}^{g_{|T|}} d_{l_g}(g)\right)$, where $m\in\{0,\cdots,N|T|\}$. Finding combinations of $(l_{g_1},\cdots,l_{g_{|T|}})$ such that $l_{g_1}+\cdots+l_{g_{|T|}}=m$ can be cast to the problem of obtaining $m$ points on $|T|$ $N+1$-sided dice (pages 23-24 in \cite{uspensky1937introduction}). The number of all possible combinations equals the coefficient of ${\epsilon}^c(t)^m$ in the polynomial $(\bm{\sum}\limits_{i=0}^N {\epsilon}^c(t)^i)^{|T|}$. By generalizing the solution of problem 13 in \cite{uspensky1937introduction}, we can get the coefficient of ${\epsilon}^c(t)^m$ is $\bm{\sum}\limits_{i=0}^{\lfloor\frac{m}{N+1}\rfloor} (-1)^i {|T| \choose i} {|T|+m-(N+1)i-1 \choose m-(N+1)i}$, where $\lfloor \cdot\rfloor$ is the floor function.
And by multinomial theorem~(Section 24.1.2 in \cite{abramowitz1964handbook}), the summation of all coefficients in $(\bm{\sum}\limits_{i=0}^N {\epsilon}^c(t)^i)^{|T|}$ equals $(N+1)^{|T|}$.\cut{ Then we use the upper bound $\bm{\sum}\limits_{i=0}^{\lfloor\frac{m}{N+1}\rfloor} (-1)^i {|T| \choose i} {|T|+m-(N+1)i-1 \choose m-(N+1)i}\le  (N+1)^{|T|}$ for presentation simplicity.}
Combining Claim 5, we have $|f^T_m|\le (N+1)^{|T|} (2^N\max\{1,||\gamma||_\infty^N\})^{|T|}$.

Note that any tree $T\in G_{{\epsilon}(t)}(z)$ is a spanning tree of the graph $\mathcal{G}({\epsilon}(t))$ where each vertex is a state $z\in \mathcal{Z}$. Then $|T|= |\mathcal{Z}|-1\le |\mathcal{Z}|$\cut{ $|T|\le {|\mathcal{Z}|(|\mathcal{Z}|-1)/2}\le |\mathcal{Z}|$}. Therefore, $|f^T_m|\le(N+1)^{|\mathcal{Z}|}2^{N|\mathcal{Z}|} C_{max} \cut{\max\{1,||\gamma||_\infty^{N|\mathcal{Z}|}\}}$.

By Cayley's formula~\cite{10027155976}, there are $|\mathcal{Z}|^{|\mathcal{Z}|-2}$ spanning trees on the complete graph with each vertex being a state $z\in \mathcal{Z}$. From (\ref{Cor1e1}) and (\ref{Cor1e4}), $b^z_m=\bm{\sum}\limits_{T\in G_{{\epsilon}(t)}(z)}f^T_m$. With $|\mathcal{Z}|^{|\mathcal{Z}|-2}\le |\mathcal{Z}|^{|\mathcal{Z}|}$, for all $m\in\{k,\cdots,h\}$, it holds that $|b^z_m|\le\bm{\sum}\limits_{T\in G_{{\epsilon}(t)}(z)}|f^T_m|\le|\mathcal{Z}|^{|\mathcal{Z}|}(N+1)^{|\mathcal{Z}|}2^{N|\mathcal{Z}|} C_{max} \cut{\max\{1,||\gamma||_\infty^{N|\mathcal{Z}|}\}}$.


Based on equalities (\ref{S4e1}), (\ref{L2c1eq1}) and (\ref{Cor1e1}), we have:
\begin{align}\label{Cor1e3}
\nonumber&b_k{\epsilon}^c(t)^k+\cdots+b_h{\epsilon}^c(t)^h=\bm{\sum}\limits_{z\in \mathcal{Z}}\sigma_{z}({\epsilon}(t))\\
&=\bm{\sum}\limits_{z\in \mathcal{Z}}b^z_k{\epsilon}^c(t)^k+\cdots+\bm{\sum}\limits_{z\in \mathcal{Z}}b^z_h{\epsilon}^c(t)^h.
\end{align}
Therefore, $|b_m|\le \bm{\sum}\limits_{z\in \mathcal{Z}}|b^z_m|\le |\mathcal{Z}|^{|\mathcal{Z}|+1}(N+1)^{|\mathcal{Z}|}2^{N|\mathcal{Z}|}\allowbreak C_{max} \cut{\max\{1,||\gamma||_\infty^{N|\mathcal{Z}|}\}}, \forall m\in\{k,\cdots,h\}$. By equality (\ref{Cor1e4}) and $|T|\le|\mathcal{Z}|$, the largest degree of $\prod\limits_{(x,y)\in E(T)}P^{{\epsilon}(t)}(x,y)$ is less than ${N|\mathcal{Z}|}$, hence $h\le N|\mathcal{Z}|$. Then $|b^z_{k+1}|+\cdots+|b^z_h|$ and $|b_{k+1}|+\cdots+|b_h|$ have less than $N|\mathcal{Z}|$ terms, which implies that $|b^z_{k+1}|+\cdots+|b^z_h|\le N|\mathcal{Z}|^{|\mathcal{Z}|+3}(N+1)^{|\mathcal{Z}|}2^{N|\mathcal{Z}|} C_{max} \cut{\max\{1,||\gamma||_\infty^{N|\mathcal{Z}|}\}}$ and $|b_{k+1}|+\cdots+|b_h|\le N|\mathcal{Z}|^{|\mathcal{Z}|+3}(N+1)^{|\mathcal{Z}|}2^{N|\mathcal{Z}|} C_{max} \cut{\max\{1,||\gamma||_\infty^{N|\mathcal{Z}|}\}}$.  \qed\end{pf}

{\bf Claim 7. $|b_k+b_{k+1}{\epsilon}^c(t) +\cdots+b_h{\epsilon}^c(t)^{h-k}|\ge \frac{1}{2}C_{min} \cut{\min\{(\frac{|\gamma|_{min}}{|\mathcal{A}|_{\infty}})^{N|\mathcal{Z}|}, 1\}}$ when $||\tilde\epsilon(t)||_\infty\le C_{min} \cut{\min\{(\frac{|\gamma|_{min}}{|\mathcal{A}|_{\infty}})^{N|\mathcal{Z}|}, 1\}}/2(N|\mathcal{Z}|^{|\mathcal{Z}|+3}(N+1)^{|\mathcal{Z}|}2^{N|\mathcal{Z}|}C_{max} \cut{\max\{1,||\gamma||_\infty^{N|\mathcal{Z}|}\}} )$.}

\begin{pf}
Now let us consider the magnitude of the coefficient of the least degree term in $\sigma_z({\epsilon}(t))$; i.e., $|b^z_k|$. Denote by $\hat f^T$  the coefficient of the least degree term in $\prod\limits_{(x,y)\in E(T)}P^{{\epsilon}(t)}(x,y)$. From equality (\ref{Cor1e1}),\cut{ $b^z_k=\bm{\sum}\limits_{T\in G_{{\epsilon}(t)}(z)} \hat f^T$ and} $\hat f^T= \prod\limits_{(x,y)\in E(T)} \hat d(x,y)$.
\cut{ and $|\mathcal{A}|_{\infty}\triangleq \max\limits_{i\in \mathcal{V}}|\mathcal{A}_i|$}%
By Claim 5, for any $T\in G_{{\epsilon}(t)}(z)$, $\hat f^T\ge \min\{\prod\limits_{(x,y)\in E(T)}(\frac{|\gamma|_{min}}{|\mathcal{A}|_{\infty}})^N\allowbreak, 1\}\ge C_{min} \cut{\min\{(\frac{|\gamma|_{min}}{|\mathcal{A}|_{\infty}})^{N|\mathcal{Z}|}, 1\}}$ because $|T|\le|\mathcal{Z}|$. \cut{Let ${T^\ast}$ be the tree with least degree among $G_{{\epsilon}(t)}(z)$. }Then $b^z_k\ge \hat f^{T}\ge C_{min} \cut{\min\{(\frac{|\gamma|_{min}}{|\mathcal{A}|_{\infty}})^{N|\mathcal{Z}|}, 1\}}$.

From equality (\ref{Cor1e3}), we have $b_k=\bm{\sum}\limits_{z\in \mathcal{Z}}b^z_k$. Since $b^z_k$ is positive, $b_k\ge b^z_k$. When\\$||\tilde\epsilon(t)||_\infty\le C_{min} \cut{\min\{(\frac{|\gamma|_{min}}{|\mathcal{A}|_{\infty}})^{N|\mathcal{Z}|}, 1\}}/2(N|\mathcal{Z}|^{|\mathcal{Z}|+3}(N+1)^{|\mathcal{Z}|}2^{N|\mathcal{Z}|}C_{max} \cut{\max\{1,||\gamma||_\infty^{N|\mathcal{Z}|}\}} )$, $|b_{k+1}{\epsilon}^c(t)+\cdots+b_h{\epsilon}^c(t)^{h-k}|\le \frac{1}{2}C_{min} \cut{\min\{(\frac{|\gamma|_{min}}{|\mathcal{A}|_{\infty}})^{N|\mathcal{Z}|}, 1\}}$ and $|b_{k}+\cdots+b_h{\epsilon}^c(t)^{h-k}|\ge \frac{1}{2}C_{min} \cut{\min\{(\frac{|\gamma|_{min}}{|\mathcal{A}|_{\infty}})^{N|\mathcal{Z}|}, 1\}}$. \qed\end{pf}

Combining Claim 6 and Claim 7, (\ref{Cor1in1}) becomes $\underline{C_\epsilon}\le C_\epsilon$.
\qed\end{pf}

\section{Case studies}\label{example}

\cut{The model in Section~\ref{model} is able to characterize the interactions between real-world defenders and attackers. For different attackers, their attack actions are different. And the corresponding defense actions, utilities are also different. In this section, we study an attacker equipped with a set zero-day attack scripts and the defender equipped with a set of platforms. Then we use the simulation based on real-world data to demonstrate its performance.}

In this section, we will evaluate the RL algorithm using two applications; i.e., the demand allocation market in~\cite{Zhu:PESGM14} and the cyber security scenario in~\cite{okhravi2014quantitative}.

\subsection{Case 1: Demand allocation market}\label{sec:case1}
In this section, we study a power market\cut{ described in Figure~\ref{fig6}. The power market} which consists of $N$ customers and a system operator. Each customer wants to allocate its demands in near future time slots and its action is subject to a price enforced by the system operator.

\subsubsection{System components}
{\bf Customers.} We consider $N$ customers $\mathcal{V}=\{1,\cdots,N\}$ and each customer $i\in \mathcal{V}$ has power demands $x_i\ge 0$ and wants to allocate its demands in one time slot within $\mathcal{A}_i=\{1,2,\cdots,|\mathcal{A}_i|\}$. The action $a^i\in \mathcal{A}_i$ is the time slot chosen by customer $i$.\cut{ Note that in this case study, $\mathcal{A}_i$ are identical for all $i\in\mathcal{V}$. At each iteration,} Each customer wants to satisfy its demands as soon as possible so it punishes late allocation. The cost function $c_i:\mathcal{A}_i\to \mathbb{R}$ is not decreasing; i.e., $c_i(a^i)\le c_i(\hat{a}^i)$ if $\hat{a}^i>a^i$.

{\bf System operator.} The system operator charges each customer some price based on demand distributions. In particular, given an action profile $s=(a^1,\cdots,a^N)$, the total demand allocated in time slot $a^i$ is $\Xi_{a^i}(s)\triangleq \bm{\sum}\limits_{j\in \mathcal{V}}\mathbf{1}_{\{a^j=a^i\}}x_j$, where $\mathbf{1}_{\{\Pi\}}$ is an indicator function: $\mathbf{1}_{\{\Pi\}}=1$ if $\Pi$ is true and $\mathbf{1}_{\{\Pi\}}=0$ if $\Pi$ is false. The system operator charges customer $i$ the price $p_a(\Xi_{a^i}(s))$.

{\bf Utility.} The utility of customer $i$ is the negative of the cost and price: $u_i(s)=-c_i(a^i)-p_a(\Xi_{a^i}(s))$.

{\bf Informational constraint.}
Each customer is unwilling to share its cost function $c_i$ and private action $a^i$ with other customers and the system operator. And the system operator does not want to disclose the pricing policy to the customers and only agrees to publicize the price value $p_a(s)$ given $s$. Therefore, each customer only knows its own utility values instead of the structure of the utility function.

\subsubsection{Evaluation}
{\bf Evaluation setup.} In this section, we use Matlab simulations to evaluate the performance of the RL algorithm. Similar to the setup in~\cite{Zhu:PESGM14}, we consider 100 customers and they have identical action sets consisting of 10 time slots. The demands of all customers are 1; i.e., $x_i=1$ for any $i\in\mathcal{V}$. The cost function for customer $i$ is set as $c_i(a^i)=\rho_i \xi_i^{a^i}$, where $\rho_i>0$ and $\xi_i> 1$. And the pricing mechanism is $p_a(\Xi_{a^i}(s))=\Xi_{a^i}(s)$.

{\bf Nash equilibrium.}
By Lemma 2.1 in~\cite{Zhu:PESGM14}, we know that the demand allocation game under the above setup is a potential game, and then a weakly acyclic game~\cite{DM-LSS:96}. Therefore the existence of pure Nash equilibrium is guaranteed. 
\begin{figure}[t]
  \centering
  \includegraphics[scale = 0.36]{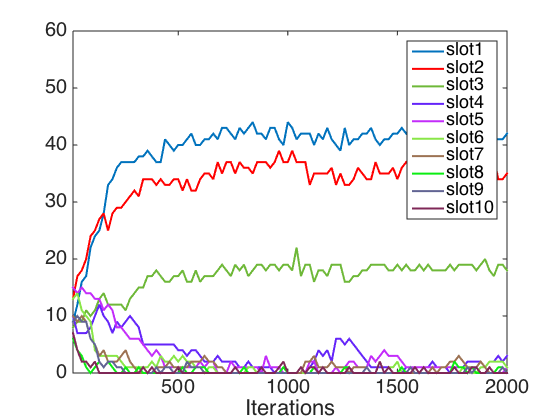}
  \caption{\cut{Each curve represents the }Temporal aggregate demands allocated at ten time slots with diminishing exploration rate $\epsilon_i(t)=\frac{1}{10}t^{-\frac{0.25}{100}}$.}\label{fig11}
\end{figure}
\begin{figure}[t]
  \centering
  \includegraphics[scale = 0.36]{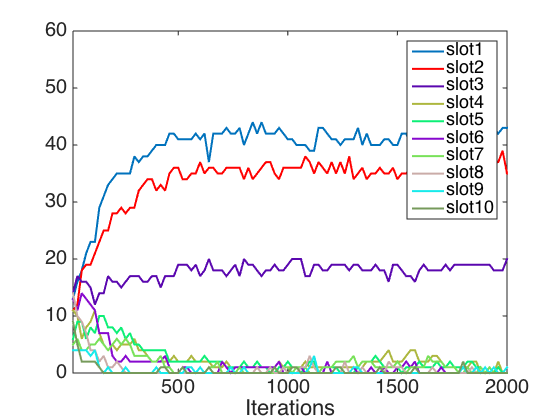}
  \caption{Temporal aggregate demands allocated at ten time slots with diminishing exploration rate $\epsilon_i(t)=\frac{1}{10}t^{-\frac{0.5}{100}}$.}\label{fig12}
\end{figure}
\begin{figure}[t]
  \centering
  \includegraphics[scale = 0.36]{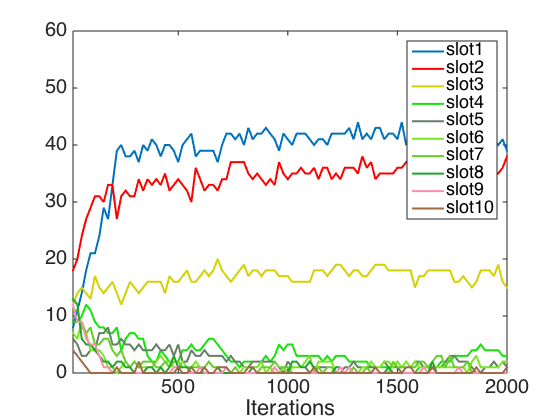}
  \caption{Temporal aggregate demands allocated at ten time slots with diminishing exploration rate $\epsilon_i(t)=\frac{1}{10}t^{-\frac{1}{100}}$.}\label{fig13}
\end{figure}

{\bf Simulation results with diminishing exploration rates.}
As discussed in Section \ref{sec:optexp}, we restrict the exploration rates to be $\mathbf{p}$-series. In particular, Figures~\ref{fig11} - \ref{fig13} show the evaluations of the RL algorithm for three cases $\epsilon_i(t)=\frac{1}{10}t^{-\frac{0.25}{100}}$, $\epsilon_i(t)=\frac{1}{10}t^{-\frac{0.5}{100}}$ and $\epsilon_i(t)=\frac{1}{10}t^{-\frac{1}{100}}$ (the optimal one), respectively. The exploration deviations are chosen as $e_i(t)=\frac{9}{10t^{2}}$ and the measurement noises are absent; i.e., $w_i(t) = 0$.
The duration of the simulation is 2,000 iterations. The simulation results confirm the convergence of the action profiles in Theorem~\ref{The4}. In addition, the convergence in Figure \ref{fig13} is fastest where the optimal exploration rates $\epsilon_i(t)=\frac{1}{10}t^{-\frac{1}{100}}$ are adopted. It is consistent with the discussion in Section~\ref{sec:optexp}.

\begin{figure}[h]
  \centering
  \includegraphics[scale = 0.36]{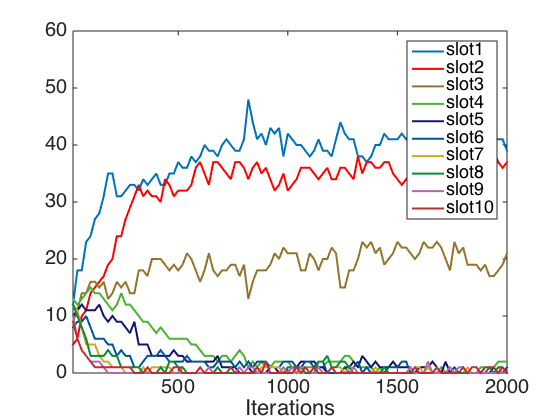}
  \caption{Temporal aggregate demands allocated at ten time slots with uniformly distributed measurement noises in the interval $[-10,10]$.}\label{fig16}
\end{figure}

\begin{figure}[h]
  \centering
  \includegraphics[scale = 0.36]{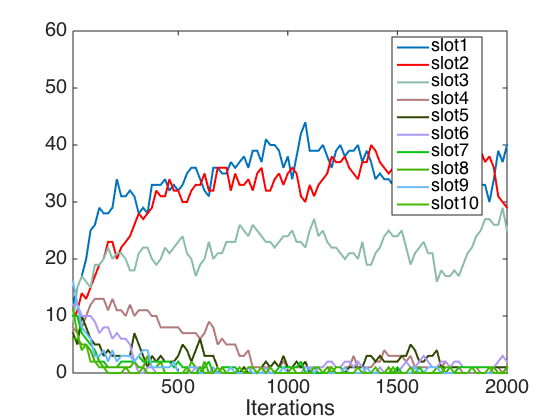}
  \caption{Temporal aggregate demands allocated at ten time slots with uniformly distributed measurement noises in the interval $[-20,20]$.}\label{fig17}
\end{figure}

\begin{figure}[h]
  \centering
  \includegraphics[scale = 0.36]{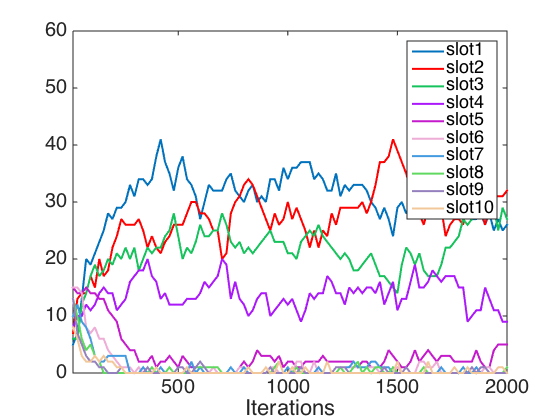}
  \caption{Temporal aggregate demands allocated at ten time slots with uniformly distributed measurement noises in the interval $[-10\ln(t),10\ln(t)]$.}\label{fig18}
\end{figure}

{\bf Simulation results with measurement errors.}
In this part, we evaluate how the RL algorithm performs when the measurement noises are present. The exploration rates are chosen as $\epsilon_i(t)=\frac{1}{10}t^{-\frac{1}{100}}$ and the exploration deviations are chosen as $e_i(t)=\frac{9}{10t^{2}}$. The measurement noises are chosen as uniformly distributed over two different intervals $[-10,10]$ and $[-20,20]$, respectively. Compared with Figure~\ref{fig13}, Figures \ref{fig16} - \ref{fig17} show that the action profiles slow down and oscillate with larger magnitudes when the noise magnitude increases. In addition, we also evaluate how the RL algorithm performs when Assumption~\ref{asu_error_o} is violated. In particular, $w_i(t)$ is uniformly distributed over the time-dependent interval $[-10\log(t),10\log(t)]$. The result shown in Figures \ref{fig18} implies that the action profiles do not converge anymore.

\begin{figure}[h]
  \centering
  \includegraphics[scale = 0.36]{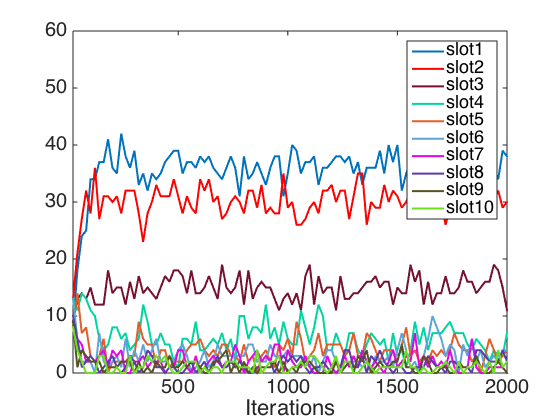}
  \caption{Temporal aggregate demands allocated at ten time slots with fixed exploration rates $\epsilon_i(t)=\frac{1}{10}^{-\frac{1}{100}}$.}\label{fig15}
\end{figure}

{\bf Matlab simulation results with fixed exploration rates.} Figures \ref{fig15} shows the evaluation of the RL algorithm with fixed exploration rates $\epsilon_i(t)=\frac{1}{10}^{-\frac{1}{100}}$. The exploration deviations are chosen as $e_i(t)=\frac{9}{10t^{2}}$ and the measurement noises are absent; i.e., $w_i(t) = 0$.  The comparison of Figures~\ref{fig13} and \ref{fig15} shows that fixed exploration rates cause larger oscillations in steady state.


\subsection{Case 2: Adaptive cyber defense scenario}

\begin{table*}\caption[Utility table]{Utility table}\label{utable}
\centering

\begin{tabular}{|c|c|c|c|c|c|c|c|c|c|c|c|}\hline

\multirow{2}{*}{Defense actions} & \multicolumn{11}{c|}{Attack actions} \\\cline{2-12}

  &(0,10) & (1,9) &(2,8) & (3,7) &(4,6) & (5,5) &(6,4) & (7,3) &(8,2) &(9,1) & (10,0) \\\hline

$d_1$(Fedora 11) &0,1.0 & 0.1,0.9 &0.2,0.8 & 0.3,0.7 &0.4,0.6 & 0.5,0.5&0.6,0.4 &0.7,0.3 &0.8,0.2 &0.9,0.1 & 1.0,0 \\\hline

$d_2$(Gentoo 9) &1.0,0 & 0.9,0.1 &0.8,0.2 & 0.7,0.3 &0.6,0.4 & 0.5,0.5&0.4,0.6 &0.3,0.7 &0.2,0.8 &0.1,0.9 & 0,1.0 \\\hline


$d_3$(CentOS 6.3) &0,1.0 & 0.1,0.9 &0.2,0.8 & 0.3,0.7 &0.4,0.6 & 0.5,0.5&0.6,0.4 &0.7,0.3 &0.8,0.2 &0.9,0.1 & 1.0,0 \\\hline
$d_4$(Debian 6) &1.0,0 & 1.0,0 &1.0,0 & 1.0,0 &1.0,0 &1.0,0 &1.0,0& 1.0,0 &1.0,0 &1.0,0 & 1.0,0\\\hline
$d_5$(FreeBSD 9) &1.0,0 & 1.0,0 &1.0,0 & 1.0,0 &1.0,0 &1.0,0 &1.0,0& 1.0,0 &1.0,0 &1.0,0 & 1.0,0 \\\hline

\end{tabular}
\end{table*}
In this section, we study a real-world cyber security scenario which consists of two players: the defender and attacker. The system is the server containing several zero-day security vulnerabilities. A zero-day attack happens once that a software/hardware vulnerability is exploited by the attacker before software the engineers develop any patch to fix the vulnerability. The attacker is equipped with a set of zero-day attack scripts denoted as $\mathcal{A}$ and the defender is equipped with a set of platforms denoted as $\mathcal{D}$. The defender uses a defensive technique called dynamic platforms~\cite{okhravi2014quantitative}, which changes the properties of the server such that it is harder for the attacker to succeed\cut{ the platform properties of the server to make the attacks more complicated}. The components in the cyber security scenario and the interactions among the components will be discussed in the following paragraphs.

\subsubsection{System components}
{\bf Defender.}
The defender has a set of different platforms; e.g., different versions of operating systems and architectures. The defender periodically restarts the server, chooses one platform from $\mathcal{D}$ and deploys it on the server each time it restarts the server. The iteration denoted in Section~\ref{general_model} is the defense period in this scenario.\cut{ We consider the defense period as the defense iteration.} The\cut{ defense action set $\mathcal{D}$ is the set of platforms and the} action $d(t)$ is the platform deployed at iteration $t$.

{\bf Attacker.}
The attacker has a set of zero-day attack scripts, where each attack script can only succeed on some platforms, but no the others. The attacker periodically chooses one of the attack scripts to attack the server. Notice that the attack period is often smaller than the defense period because the defender cannot restart the server too frequently due to the resource consumption of restarting the server. In fact, the defense period is usually a multiple of the attack period. The attack action at iteration $t$, denoted as $a(t)$, is a subset of the attack scripts and the attack action set $\mathcal{A}$ includes all possible subsets. The order of choosing the attack scripts in one iteration does not matter. \cut{The attack action set $\mathcal{A}$ is the set of combinations of attack scripts determined at the beginning of $t$ and the action $a(t)$ of defense iteration $t$ is a fixed number of predetermined attack scripts (one combination).}

{\bf Server.}
Once an attack action succeeds, the attacker can control the server for a certain amount of time.\cut{ until it launches next attack script (to achieve different security goal) or the server restarts.} And every time the server restarts, the defender takes over the control of the server. And here we assume the time consumed by the attack scripts to succeed and the time consumed by restarting the server are negligible compared with the length of an iteration.

{\bf Utility.} The goal of both the attacker and defender is to gain longer control time of the server. The utility of the defender $u_d(d(t),a(t))$ is the fraction of the time controlled by the defender during iteration $t$ and the utility of the attacker $u_a(d(t),a(t))$ is the fraction of the  time controlled by the attacker. Notice that $u_d(d(t),a(t))+u_a(d(t),a(t))=1$.

{\bf Informational constraint.}
The attacker can observe when the server restarts, so it knows the iterations, but it does not know the defender's action set $\mathcal{D}$ and which platform is deployed. The defender does not know the attacker's action set $\mathcal{A}$ and the specific attack scripts chosen by the attacker. At the end of each defense period, both the defender and attacker can measure how much time they control the server. Therefore, each player only knows its own utility values instead of the structure of the utility function.

\subsubsection{Evaluation}
{\bf Evaluation setup.} In this section, we use Matlab simulations to evaluate the performance of our algorithm based on real-world platform settings, attack scripts and server control data~\citep{okhravi2012creating,okhravi2014quantitative}. The total number of defense actions is 5; i.e., the defender has five different platforms: Fedora 11 on x86, Gentoo 9 on x86, Debian 6 on x86, FreeBSD 9 on x86, and CentOS 6.3 on x86. The attacker has two zero-day attack scripts: TCP MAXSEG exploit, and Socket Pairs exploit. The defense period is set to be ten times as large as the attack period; i.e., during one iteration, the attacker launches 10 attack scripts. Since the time consumed by the attack scripts to succeed is negligible, one attack script enables the attacker control $\frac{1}{10}$ of the iteration if it succeeds.\cut{ Then the order of choosing the attack scripts in one defense iteration does not matter and the total number of attack actions is 11}  The total number of attack actions is 11; i.e., $a_1=(0,10),a_2=(1,9),\cdots,a_{11}=(10,0)$, where $a(t)=(0,10)$ means the attacker launches 0 TCP MAXSEG exploit and 10 Socket Pairs exploits at iteration $t$.

{\bf Real-world utility values.} Based on the evaluation setup and the real-world attack scripts, we first replay different attack actions on different platforms to get the utility table for the defender and the attacker. The results are shown in Table~\ref{utable}, where the defender is the row player and the attacker is the column player. In each cell, the first number represents the utility value to the defender, and the second number represents the utility value to the attacker.

{\bf Nash equilibrium.}
By Proposition 1 in~\cite{takahashi2002pure}, we know any 2-player finite game and its any sub-game (any game constructed by restricting the set of actions to a subset of the set of actions in the original game) has at least one pure Nash equilibrium is a weakly acyclic game. From Table~\ref{utable}, we can see any sub-game has at least one pure Nash equilibrium. Now we want to calculate the pure Nash equilibrium (equilibria). From Table~\ref{utable}, we can see if the defense strategy is $d_4$ (deploying Debian 6) or $d_5$ (deploying FreeBSD 9), then the utility of the defender is 1 (the utility of the attacker is 0) not matter what action the attacker uses. From Definition~\ref{def1_NE} and Remark~\ref{rmk1}, we know the combinations of any attacker action and defense action $d_4$ or $d_5$ are pure Nash equilibria.
%
%
%
\begin{figure}[h]
  \centering
  \includegraphics[scale = 0.35]{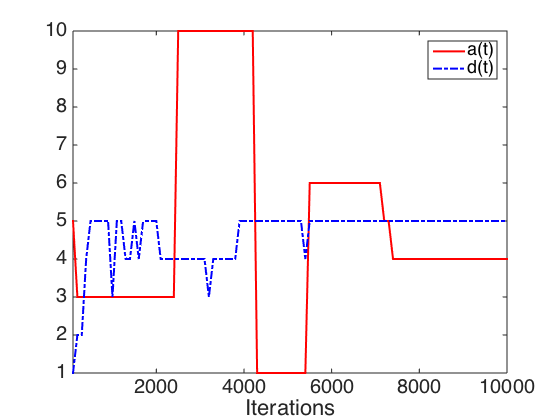}
  \caption{Trajectories of attack $a(t)$ and defense $d(t)$ with diminishing exploration rate $\epsilon_d(t)= \epsilon_a(t)=\frac{1}{11t^{1/2}}$.}\label{fig4}
\end{figure}


{\bf Simulation results with optimal exploration rates.} Based on Table~\ref{utable}, we simulate the interactions of the defender and attacker in Matlab.
We choose the exploration rates $\epsilon_d(t)= \epsilon_a(t)=\frac{1}{11t^{1/2}}$., the optimal one among $\mathbf{p}$-series as discussed in Section \ref{sec:optexp}. The exploration deviations are chosen as $e_d(t)=\frac{1}{110t^{2}}$ and $e_a(t)=\frac{1}{110t^{2}}$. We assume that the measurement noises are absent; i.e., $w_d(t)=w_a(t) = 0$.
The duration of each simulation (from the attack begins till the attack ends) is 10,000 iterations and we repeat 100 identical simulations. Figure~\ref{fig4} shows the trajectories of the defense and attack actions in one certain simulation. And for each simulation, we record the defense action at each iteration. Then at each iteration $t$, we have 100 chosen defense actions and we use the number of each defense action over 100 as the probability of choosing such defense action at $t$. The result in Figure~\ref{fig4} suggests that the defense action converges to the set $\{d_4,d_5\}$. 
Notice that the combinations of any attacker action and defense action $d_4$ or $d_5$ are pure Nash equilibria. Then the simulation results confirm that the convergence of the action profiles to the set of pure Nash equilibria.

\section{Conclusion}
This paper investigates a class of multi-player discrete games where each player aims to maximize its own utility function with limited information about the game of interest. We propose the RL algorithm which converges to the set of action profiles which have maximal stochastic potential with probability one. The convergence rate of the proposed algorithm is analytically quantified. Moreover, the performance of the algorithm is verified by two case studies in the smart grid and cybersecurity.

\bibliographystyle{elsarticle-harv}
\bibliography{REF_TAAS_TIFS_TDSC_AUTOMATICA_201710,REF_TAAS_TIFS_TDSC_AUTOMATICA_201610,MZ,Zhu}           
\cut{{\bf Differences between the mathematical model in this paper and the one in the S\&P paper:}

{\bf Allocator:} In this paper, the buffers are randomly located on the heap by following certain determined probability distribution; i.e., $\vec{SA}$ is a random vector. In the S\&P paper, the buffers are randomly located on the heap by following some probability distributions conditional on the defense actions.

{\bf Attacker:} In this paper, the attacker has a diminishing exploration rate $\epsilon_a(t)$, and with probability $1-\epsilon_a(t)$, the attacker sticks to the previous action $a(t-1)$; with probability $\epsilon_a (t)$, it switches to a new one via following its algorithm $ALG_a$. In the S\&P paper, the attacker has a set of predetermined probability distributions $\{DA_d\}$ associated with defense actions. And the attacker updates attack actions simultaneously with the defender: firstly, the attacker precisely observes the defender's current action $d(t)$ at the beginning of each defense iteration $t$; secondly, given an observation $d(t)$, the attacker chooses its action $a(t)$ according to the corresponding probability distribution $DA_{d(t)}$.}


\end{document}